\documentclass[12pt]{article}
\usepackage[mathscr]{eucal}
\usepackage{amssymb}
\usepackage{latexsym}
\usepackage{amsmath}
\newtheorem{definition}{Definition}[section]
\newtheorem{lemma}[definition]{Lemma}
\newtheorem{theorem}[definition]{Theorem}
\newtheorem{proposition}[definition]{Proposition}
\newtheorem{corollary}[definition]{Corollary}

\begin{document}
\title{Cuntz-Krieger-Pimsner Algebras \\ Associated with \\ Amalgamated Free
Product Groups}
\author{Rui OKAYASU \\ Department of Mathematics, Kyoto University, \\
Sakyo-ku, Kyoto, 606-8502, Japan \\ {\tt e-mail:rui@kusm.kyoto-u.ac.jp}}
\date{}

\maketitle

\begin{abstract}
{\footnotesize We give a construction of a nuclear
$C^\ast$-algebra associated with an amalgamated free
product of groups, generalizing Spielberg's construction of a certain
Cuntz-Krieger algebra associated with a finitely generated free product of cyclic
groups. Our nuclear $C^\ast$-algebras can
be identified with certain Cuntz-Krieger-Pimsner algebras. We will also
show that our algebras can be obtained by the crossed product
construction of the canonical actions on the hyperbolic boundaries, which
proves a special case of Adams' result about amenability of the boundary action  for 
hyperbolic groups. We will also give an explicit formula of the
$K$-groups of our algebras. Finally we will investigate the
relationship between the KMS states of the generalized gauge actions on our
$C^\ast$ algebras and random walks on the groups.}
\end{abstract}

\section{Introduction}

In ~\cite{cho}, Choi proved that the reduced group $C^\ast$-algebra
$C_r^\ast\left(\mathbb{Z}_2\ast\mathbb{Z}_3\right)$ of the free
product of cyclic groups $\mathbb{Z}_2$ and $\mathbb{Z}_3$ is embedded in $\mathcal{O}_2$. Consequently, this shows that 
$C_r^\ast\left(\mathbb{Z}_2\ast\mathbb{Z}_3\right)$ is a non-nuclear
exact $C^\ast$-algebra, (see S. Wassermann ~\cite{was} for a good introduction to exact $C^\ast$-algebras). Spielberg 
generalized it to finitely generated free products of cyclic
groups in ~\cite{spi}. Namely, he constructed a certain action on a compact space and proved
that some Cuntz-Krieger algebras (see ~\cite{ck1}) can be obtained by the crossed
product construction for the action. For a related topic, see
W. Szyma\'nski and S. Zhang's work ~\cite{sz}.

More generally, the above mentioned compact space coincides with
Gromov's notion of the boundaries of hyperbolic groups (e.g. see
~\cite{gph}). In ~\cite{ada}, Adams
proved that the action of any discrete hyperbolic group $\Gamma$ on
the hyperbolic boundary $\partial\Gamma$
is amenable in the sense of Anantharaman-Delaroche ~\cite{ana}. It follows from ~\cite{ana}
that the corresponding crossed product
$C(\partial\Gamma)\rtimes_r\Gamma$ is nuclear, and this
implies that $C_r^\ast(\Gamma)$ is an exact $C^\ast$-algebra.

Although we know that $C(\partial\Gamma)\rtimes_r\Gamma$ is nuclear for a general 
discrete hyperbolic group $\Gamma$ as mentioned above, there are only few things known about 
this $C^\ast$-algebra. So one of our purposes is to generalize Spielberg's construction to
some finitely generated amalgamated free product $\Gamma$ and to
give detailed description of the algebra $C(\partial\Gamma)\rtimes_r\Gamma$. More precisely, 
let $I$ be
a finite index set and $G_i$ be a group containing a copy of a finite group $H$ as a subgroup for $i\in I$. We always assume that
each $G_i$ is either a finite group or $\mathbb{Z}\times
H$. Let $\Gamma=\ast_H G_i$ be the amalgamated free product group. We will construct a nuclear $C^\ast$-algebra
$\mathcal{O}_{\Gamma}$ associated with $\Gamma$ by mimicking the construction 
for Cuntz-Krieger algebras with respect to the full Fock space in M. Enomoto, M. Fujii 
and Y. Watatani ~\cite{efw} and D. E. Evans ~\cite{eva}. This
generalizes Spielberg's construction. 

First we show that $\mathcal{O}_{\Gamma}$ has a certain universal
property as in the case of the Cuntz-Krieger algebras, which allows
several descriptions of $\mathcal{O}_{\Gamma}$. For example, it turns
out that $\mathcal{O}_{\Gamma}$ is a Cuntz-Krieger-Pimsner algebra,
introduced by Pimsner in ~\cite{pim2} and studied by several authors, e.g. T. Kajiwara, C. Pinzari and Y. Watatani ~\cite{kpw}. We will also show
that $\mathcal{O}_{\Gamma}$ can be obtained by the crossed product
construction. Namely, we will introduce a boundary space $\Omega$ with 
a natural $\Gamma$-action, which coincides with the boundary of the
associated tree (see ~\cite{ser}, ~\cite{woe}). Then we will prove that $C(\Omega)\rtimes_r\Gamma$ is isomorphic to $\mathcal{O}_{\Gamma}$. 
Since the hyperbolic boundary $\partial\Gamma$ coincides 
with $\Omega$ and the two actions of $\Gamma$ on $\partial\Gamma$ and $\Omega$ are 
conjugate, $\mathcal{O}_{\Gamma}$ is also isomorphic to 
$C(\partial\Gamma)\rtimes_r\Gamma$, and depends only on the group
structure of $\Gamma$. As a consequence, we give a proof to Adams' theorem in 
this special case.

Next, we will consider the $K$-groups of $\mathcal{O}_{\Gamma}$. In ~\cite{pim}, Pimsner 
gave a certain exact sequence of $KK$-groups
of the crossed product by groups acting on trees. However, it is not a 
trivial task to apply Pimsner's exact sequence to
$C(\partial\Gamma)\rtimes_r\Gamma$ and obtain its $K$-groups. We will give explicit
formulae of the $K$-groups of $\mathcal{O}_{\Gamma}$ following the
method used for the Cuntz-Krieger algebras instead of using
$C(\partial\Gamma)\rtimes_r\Gamma$. We can compute the $K$-groups of 
$C(\partial\Gamma)\rtimes_r\Gamma$ for concrete examples. They are completely determined by 
the
representation theory of $H$ and the actions of $H$ on $G_i/H$ (the space of right cosets) by left multiplication.

Finally we will prove that KMS states on $\mathcal{O}_{\Gamma}$ for
generalized gauge actions arise from harmonic measures on the Poisson
boundary with respect to random walks on the discrete group
$\Gamma$. Consequently, for special cases, we can determine easily the type of factor
$\mathcal{O}_{\Gamma}''$ for the corresponding unique KMS state of the
gauge action by essentially the same arguments in M. Enomoto, M. Fujii 
and Y. Watatani ~\cite{efw2}, which generalized J. Ramagge and
G. Robertson's result ~\cite{rr}.

\medskip

{\bf Acknowledgment.} The author gives special thanks to Professor Masaki Izumi for various 
comments and many important suggestions.

\section{Preliminaries}
In this section, we collect basic facts used in the present article.
We begin by reviewing the Cuntz-Krieger-Pimsner algebras in
~\cite{pim2}. Let $A$ be a $C^\ast$-algebra and $X$ be a Hilbert
bimodule over $A$, which means that $X$ is a right Hilbert $A$-module
with an injective $\ast$-homomorphism of $A$ to $\mathcal{L}(X)$, where $\mathcal{L}(X)$ is 
the $C^\ast$-algebra of all adjointable $A$-linear operators on $X$. We assume
that $X$ is full, that is, $\{\langle x, y\rangle_A \mid x,
y\in X\}$ generates $A$ as a $C^\ast$-algebra, where $\langle\cdot
, \cdot\rangle_A$ is the $A$-valued inner product on $X$. We further
assume that $X$ has a finite basis $\{u_1, \ldots , u_n\}$, which means
that $x=\sum_{i=1}^nu_i\langle u_i, x\rangle_A$ for any $x\in X$. We 
fix a basis $\{u_1, \ldots , u_n\}$ of $X$. Let
$\mathcal{F}(X)=A\oplus\bigoplus_{n\geq 1}X^{(n)}$ be the full Fock
space over $X$, where $X^{(n)}$ is the $n$-fold tensor product
$X\otimes_AX\otimes_A\cdots\otimes_A X$. Note that $\mathcal{F}(X)$ is naturally equipped with Hilbert $A$-bimodule structure. For each $x\in X$, the operator $T_x : \mathcal{F}(X)\to\mathcal{F}(X)
$ is defined by
\begin{eqnarray*}
T_x(x_1\otimes\cdots\otimes x_n) &=& x\otimes x_1\otimes\cdots\otimes
x_n, \\
T_x(a) &=& xa,
\end{eqnarray*}
for $x, x_1, \ldots , x_n\in X$ and $a\in A$. Note that
$T_x\in\mathcal{L}(\mathcal{F}(X))$ satisfies the following relations
\begin{eqnarray*}
T_x^\ast T_y &=& \langle x, y\rangle_A, \qquad x, y\in X, \\
aT_xb &=& T_{axb}, \qquad x\in X, a, b\in A.
\end{eqnarray*}

Let $\pi$ be the quotient map of $\mathcal{L}(\mathcal{F}(X))$ onto
$\mathcal{L}(\mathcal{F}(X))/\mathcal{K}(\mathcal{F}(X))$ where
$\mathcal{K}(\mathcal{F}(X))$ is the $C^\ast$-algebra of all compact
operators of $\mathcal{L}(\mathcal{F}(X))$. We denote $S_x=\pi(T_x)$ for $x\in
X$. Then we define the Cuntz-Krieger-Pimsner algebra $\mathcal{O}_X$
to be 
$$\mathcal{O}_X=C^\ast(S_x \mid x\in X).$$
Since $X$ is full, a copy of $A$ acting by left multiplication
on $\mathcal{F}(X)$ is contained in $\mathcal{O}_X$. Furthermore we
have the relation
$$\sum_{i=1}^nS_{u_i}S_{u_i}^\ast=1. \eqno (\dag)$$

On the other hand, $\mathcal{O}_X$ is characterized as the universal $C^\ast$-algebra generated by 
$A$ and $S_x$, satisfying the above relations ~\cite[Theorem
3.12]{pim2}. More precisely, we have

\begin{theorem}[{~\cite[Theorem 3.12]{pim2}}]
Let $X$ be a full Hilbert $A$-bimodule and $\mathcal{O}_X$ be the
corresponding Cuntz-Krieger-Pimsner algebra. Suppose that $\{u_1,
\cdots , u_n\}$ is a finite basis for $X$. If $B$ is a $C^\ast$-algebra generated by 
$\{s_x\}_{x\in X}$ satisfying 
\begin{eqnarray*}
s_x+s_y &=& s_{x+y}, \qquad x\in X, \\
as_xb &=& s_{axb}, \qquad x\in X, a, b\in A, \\
s_x^\ast s_y &=& \langle x, y\rangle_A, \qquad x, y\in X, \\
\sum_{i=1}^ns_{u_i}s_{u_i}^\ast &=& 1. \\
\end{eqnarray*}
Then there exists a unique surjective $\ast$-homomorphism from
$\mathcal{O}_X$ onto $C^\ast(s_x)$ that maps $S_x$ to $s_x$.
\end{theorem}

Next we recall the notion of amenability for discrete
$C^\ast$-dynamical systems introduced by C. Anantharaman-Delaroche in ~\cite{ana}. Let $(A, 
G, \alpha)$ be a
$C^\ast$-dynamical system, where $A$ is a $C^\ast$-algebra, $G$ is a
group and $\alpha$ is an action of $G$ on $A$. An $A$-valued function $h$ on $G$ is said
to be of $positive$ $type$ if the matrix
$[\alpha_{s_i}(h(s_i^{-1}s_j))]\in M_n(A)$ is positive for any $s_1,
\ldots , s_n\in G$. We assume that $G$ is discrete. Then $\alpha$ is
said to be $amenable$ if there exists a net $(h_i)_{i\in I}\subset C_c(G, Z(A''))$ of functions of positive type such that 
$$
\left\{\begin{array}{ll}
 h_i(e)\leq 1 & \mbox{for}\ i\in I, \\
 \displaystyle\lim_ih_i(s)=1 & \mbox{for}\ s\in G,
\end{array}\right. $$
where the limit is taken in the $\sigma$-weak topology in the
enveloping von Neumann algebra $A''$ of $A$. We remark that this is one of
several equivalent conditions given in ~\cite[Th\'{e}or\`{e}me
3.3]{ana}. We will use the following theorems without a proof.

\begin{theorem}[{~\cite[Th\'{e}or\`{e}me 4.5]{ana}}]
Let $(A, G, \alpha)$ be a $C^\ast$-dynamical system such that $A$ is
nuclear and $G$ is discrete. Then the following are equivalent:

1) The full $C^\ast$-crossed product $A\rtimes_{\alpha}G$ is nuclear;

2) The reduced $C^\ast$-crossed product $A\rtimes_{\alpha r}G$ is nuclear;

3) The $W^\ast$-crossed product $A''\rtimes_{\alpha w}G$ is injective;

4) The action $\alpha$ of $G$ on $A$ is amenable.
\end{theorem}

\begin{theorem}[{~\cite[Th\'{e}or\`{e}me 4.8]{ana}}]
Let $(A, G, \alpha)$ be an amenable $C^\ast$-dynamical system such that $G$ is
discrete. Then the natural quotient map from $A\rtimes_{\alpha}G$ onto
$A\rtimes_{\alpha r}G$ is an isomorphism.
\end{theorem}

Finally, we review the notion of the strong boundary
actions in ~\cite{spi2}. Let $\Gamma$ be a discrete group acting by 
homeomorphisms on
a compact Hausdorff space $\Omega$. Suppose that $\Omega$ has at
least three points. The action of $\Gamma$ on $\Omega$ is said to be a
$strong$ $boundary$ $action$ if for every pair $U, V$ of non-empty open 
subsets of $\Omega$ there exists $\gamma\in\Gamma$ such that $\gamma
U^c\subset V$. The action of $\Gamma$ on $\Omega$ is said to
be $topologically$ $free$ in the sense of ~\cite{as} if the fixed point set of each 
non-trivial element
of $\Gamma$ has empty interior.

\begin{theorem}[{~\cite[Theorem 5]{spi2}}]
Let $(\Omega, \Gamma)$ be a strong boundary action where $\Omega$ is
compact. We further assume that the action is topologically
free. Then $C(\Omega)\rtimes_r\Gamma$ is purely infinite and simple.
\end{theorem}

\section{A motivating example}

Before introducing our algebras, we present a simple case of
Spielberg's construction for $\mathbb{F}_2=\mathbb{Z}\ast\mathbb{Z}$ with generators $a$ and
$b$ as a motivating example. See also ~\cite{rs}. The Cayley graph of 
$\mathbb{F}_2$ is a homogeneous 
tree of degree 4. The boundary $\Omega$ of the tree in the sense of ~\cite{fre} (see
also ~\cite{fur}) can be thought of as the set of all infinite reduced words 
$\omega=x_1x_2x_3\cdots$, where $x_i\in S=\{a, b, a^{-1},
b^{-1}\}$. Note that $\Omega$ is compact in the relative topology of
the product topology of $\prod_{\mathbb{N}}S$. In an appendix, several facts
about trees are collected for the convenience of the reader, (see also ~\cite{fig}). Left multiplication of $\mathbb{F}_2$ on $\Omega$ induces 
an action of $\mathbb{F}_2$ on $C(\Omega)$. For $x\in\mathbb{F}_2$, let
$\Omega(x)$ be the set of infinite words beginning with $x$. We identify the implementing 
unitaries in the full crossed product $C(\Omega)\rtimes\mathbb{F}_2$ with elements 
of $\mathbb{F}_2$. Let $p_x$ denote the projection defined by the
characteristic function $\chi_{\Omega(x)}\in C(\Omega)$. Note that for
  each $x\in S$,
$$p_x+xp_{x^{-1}}x^{-1}=1,$$
$$p_a+p_{a^{-1}}+p_b+p_{b^{-1}}=1,$$
hold. For $x\in S$, let $S_x\in C(\Omega)\rtimes\mathbb{F}_2$ be
a partial isometry
$$S_x=x(1-p_{x^{-1}}).$$
Then we have
$$S_x^\ast S_y=x^{-1}p_xp_yy=\delta_{x, y}S_x^\ast S_x=\delta_{x, y}(1-p_{x^{-1}}),$$
$$S_xS_x^\ast=x(1-p_{x^{-1}})x^{-1}=p_x,$$
$$S_x^\ast S_x=1-p_{x^{-1}}=\sum_{y\ne x^{-1}}S_yS_y^\ast.$$
These relations show that the partial isometries $S_x$ generate the
Cuntz-Krieger algebra $\mathcal{O}_A$ ~\cite{ck1}, where
$$A=\left(\begin{array}{rrrr}
1 & 0 & 1 & 1 \\
0 & 1 & 1 & 1 \\
1 & 1 & 1 & 0 \\
1 & 1 & 0 & 1 
\end{array}\right).$$
On the other hand, we can recover the generators of $C(\Omega)\rtimes\mathbb{F}_2$ by setting
$$x=S_x+S_{x^{-1}}^\ast \quad\mbox{and}\ \quad p_x=S_xS_x^\ast.$$
Hence we have $C(\Omega)\rtimes\mathbb{F}_2\simeq\mathcal{O}_A$.

Next we recall the Fock space realization of the Cuntz-Krieger algebras, (e.g. see 
~\cite{eva}, ~\cite{efw}). Let $\{e_a, e_b, e_{a^{-1}}, e_{b^{-1}}\}$
be a basis of $\mathbb{C}^4$. We define the Fock space associated with 
the matrix $A$ by
$$\mathcal{F}_A=\mathbb{C}e_0\oplus\bigoplus_{n\geq
  1}\left(\overline{\rm span}\{e_{x_1}\otimes\cdots\otimes e_{x_n}\mid 
  A(x_i, x_{i+1})=1\}\right),$$
where $e_0$ is the vacuum vector. For any $x \in S$, let $T_x$ be the creation
  operator on $\mathcal{F}$, given by
\begin{eqnarray*}
T_xe_0 &=& e_x, \\
T_x(e_{x_1}\otimes\cdots\otimes e_{x_n}) &=& \left\{\begin{array}{ll}
                                 e_x\otimes e_{x_1}\otimes\cdots\otimes e_{x_n} & \mbox{if}\
                                 A(x, x_1)=1, \\
                                 0 & \mbox{otherwise}.\
                          \end{array}\right. 
\end{eqnarray*}
Let $p_0$ be the rank one projection on the vacuum vector $e_0$. Note
that we have
$$T_aT_a^\ast+T_bT_b^\ast+T_{a^{-1}}T_{a^{-1}}^\ast+T_{b^{-1}}T_{b^{-1}}^\ast+p_0=1.$$
If $\pi$ is the quotient map of $\mathcal{B}(\mathcal{F})$ onto the Calkin algebra 
$\mathcal{Q}(\mathcal{F})$,
then the $C^\ast$-algebra generated by the partial isometries
$\{\pi(T_a), \pi(T_b), \pi(T_{a^{-1}}), \pi( T_{b^{-1}})\}$ is isomorphic to the
Cuntz-Krieger algebra $\mathcal{O}_A$.

Now we look at this construction from another point of view. We can perform the following 
natural identification:
$$\mathcal{F}\ni \left. \begin{array}{ccc}
    e_0 & \longleftrightarrow & \delta_e \\
    e_{x_1}\otimes\cdots\otimes e_{x_n} & \longleftrightarrow &
    \delta_{x_1\cdots x_n}
    \end{array}\right.
\in l^2(\mathbb{F}_2).$$
Under this identification, the creation
operator $T_x$ on $l^2(\mathbb{F}_2)$ can be expressed as
 \begin{eqnarray*}
T_x\delta_e &=& \lambda_x\delta_e, \\
T_x\delta_{x_1\cdots x_n} &=& \left\{\begin{array}{ll}
                                 \lambda_x\delta_{x_1\cdots x_n} &
                                 \mbox{if}\ x\ne x_1^{-1}, \\
                                 0 & \mbox{otherwise}.\
                          \end{array}\right. 
\end{eqnarray*}
where $\lambda$ is the left regular representation of $\mathbb{F}_2$. 

For a reduced word $x_1\cdots x_n\in\mathbb{F}_2$, we define the
length function $|\cdot|$ on $\mathbb{F}_2$ by $|x_1\cdots x_n|=n$. Let $p_n$ be 
the projection onto the closed linear span of $\{\delta_{\gamma}\in l^2(\mathbb{F}_2)\mid
|\gamma|=n\}$. Then we can express $T_x$ for $x\in S$ by
$$T_x=\sum_{n\geq 0}p_{n+1}\lambda_xp_n.$$

Note that this expression makes sense for every finitely generated group. In the next 
section, we generalize this construction to amalgamated free product groups. 

\section{Construction of a nuclear $C^\ast$-algebra $\mathcal{O}_{\Gamma}$}

In what follows, we always assume that $I$ is a finite index set and
$G_i$ is a group containing a copy of a finite group $H$ as a
subgroup for $i\in I$. Moreover, we assume that each $G_i$
is either a finite group or $\mathbb{Z}\times H$. We set $I_0=\{i\in I 
\mid |G_i|<\infty\}$. Let
$\Gamma=\ast_H G_i$ be the amalgamated free product. 

First we introduce a ``length 
function'' $|\cdot|$ on each $G_i$. If
$i\in I_0$, we set $|g|=1$ for any $g\in G_i\setminus
H$ and $|h|=0$ for any $h\in H$. If $i\in I\setminus I_0$ we set $|(a_i^n, h)|=|n|$ for any
$(a_i^n, h)\in G_i=\mathbb{Z}\times H$ where $a_i$ is a generator of
$\mathbb{Z}$. Now we extend the length function to $\Gamma$. Let
$\Omega_i$ be a set of left representatives of $G_i/H$ with
$e\in\Omega_i$. If $\gamma\in\Gamma$ is written uniquely as $g_1\cdots g_n h$, where
$g_1\in\Omega_{i_1}, \ldots , g_n\in\Omega_{i_n}$ with
$i_1 \ne i_2, \ldots ,i_{n-1}\ne i_n$(we write simply 
$i_1 \ne
\cdots \ne i_n$), then we define 
$$|\gamma|=\sum_{k=1}^n|g_k|.$$ 

Let $p_n$ be the projection of
$l^2\left(\Gamma\right)$ onto $l^2\left(\Gamma_n\right)$ for each $n$, 
where $\Gamma_n=\{ \,\gamma\in\Gamma \mid |\gamma|=n \, \}$.
We define partial isometries and unitary operators on $l^2 \left(
 \Gamma \right)$ by
$$ \left\{\begin{array}{ll}
   T_g=\sum_{n\geq0} p_{n+1} \lambda_g p_n & \mbox{if}\  g\in\bigcup_{i\in I}G_i
\setminus H, \\
   V_h=\lambda_h & \mbox{if}\ h\in H,
\end{array}\right. $$
where $\lambda$ is the left regular representation of $\Gamma$. Let $\pi$ be the quotient 
map of $\mathcal{B}(l^2(\Gamma))$ onto
$\mathcal{B}(l^2(\Gamma))/\mathcal{K}(l^2(\Gamma))$, where
$\mathcal{B}(l^2(\Gamma))$ is the $C^\ast$-algebra of all bounded linear
operators on $l^2(\Gamma)$ and $\mathcal{K}(l^2(\Gamma))$ is the
$C^\ast$-subalgebra of all compact operators of
$\mathcal{B}(l^2(\Gamma))$. We set $\pi(T_g)=S_g$ and $\pi(V_h)=U_h$. For $\gamma\in\Gamma$, 
we define $S_{\gamma}$ by
$$S_{\gamma}=S_{g_1} \cdots S_{g_n},$$
where $\gamma=g_1 \cdots g_n$ for some $g_1\in G_{i_1}\setminus 
H, \ldots
, g_n\in G_{i_n}\setminus H$ with $i_1\ne\cdots\ne i_n$. Note that 
$S_{\gamma}$ does not depend on the expression $\gamma=g_1\cdots g_n$. We denote
the initial projections of $S_{\gamma}$ by $Q_{\gamma}=S_{\gamma}^{\ast}\cdot S_{\gamma}$
 and the range projections by $P_{\gamma}=S_{\gamma}\cdot
S_{\gamma}^{\ast}$ for $\gamma\in\Gamma$. 

We collect several relations, which the family $\{ \, S_g , U_h \mid
g\in\bigcup_{i\in I}G_i\setminus H, h\in H \, \}$ satisfies.

For $g, g'\in\bigcup_i G_i\setminus H$ with $|g|=|g'|=1$ and $h\in H,$
$$ S_{gh}=S_g \cdot U_h, \qquad S_{hg}=U_h \cdot S_g, \eqno(1)$$
$$P_g \cdot P_{g'}=\left\{\begin{array}{ll}
                         P_g=P_{g'} & \mbox{if}\ \, gH=g'H, \\
                         0 & \mbox{if}\ \, gH\ne g'H.
                          \end{array}\right. \eqno(2)$$
Moreover, if $g\in G_i\setminus H$ and $i\in I_0$, then
$$Q_g=\sum_{\stackrel{\mbox{$\scriptstyle j\in I_0$}}
                    {\mbox{$\scriptstyle j\ne i$}}}
                    \sum_{g'\in\Omega_j\setminus\{e\}}P_{g'}
                   +\sum_{j\in I \setminus I_0 }P_{a_j}+P_{a_j^{-1}}, 
\eqno(3)$$
and if $g=a_i^{\pm 1}$ and $i\in I\setminus I_0$, then
$$Q_{a_i^{\pm 1}}=\sum_{j\in
  I_0}\sum_{g'\in\Omega_j\setminus\{e\}}P_{g'}+
\sum_{\stackrel{\mbox{$\scriptstyle j\in I\setminus I_0$}}
               {\mbox{$\scriptstyle j\ne i$}}}
\left(P_{a_j}+P_{a_j^{-1}}\right)+P_{a_i^{\pm 1}}. \eqno(3)'$$
Finally,
$$1=\sum_{i\in I_0}
    \sum_{g\in\Omega_i\setminus\{e\}}P_{g}+
    \sum_{i\in I\setminus I_0}\left(P_{a_i}+P_{a_i^{-1}}\right). \eqno(4)
$$

Indeed, (1) follows from the relations $T_{gh}=T_gV_h$
and $T_{hg}=V_hT_g$. From the definition, we have $T_{g'}^\ast T_g=\sum_{n\geq
  0}p_{n}\lambda_{g'}^\ast p_{n+1}\lambda_gp_n$. This can be non-zero 
  if and only if $|{g'}^{-1}g|=0$, i.e.  ${g'}^{-1}g\in H$. We have (2) immediately. 
The relation 
$$1=\sum_{i\in
  I_0}\sum_{g\in\Omega_i}T_gT_g^\ast+\sum_{i\in I\setminus 
  I_0}\left(T_{a_i}T_{a_i}^\ast+T_{a_i^{-1}}T_{a_i^{-1}}^\ast\right)+p_0,$$
implies (4). By multiplying $S_g^\ast$ on the left and $S_g$ on the right of equation (4) respectively, we obtain (3).

Moreover, the following condition holds:
Let
$P_i=\sum_{g\in\Omega_i}P_g$ for $i\in I_0$, and 
$P_i=P_{a_i}+P_{a_i^{-1}}$ for $i\in I\setminus I_0$. For every 
$i\in I$, we have
$$C^\ast(H)\simeq C^\ast\left(P_iU_hP_i \mid h\in
  H\right). \eqno(5)$$
Indeed, since the unitary representation $P'_iV_hP'_i$ contains the left
  regular representation of $H$ with infinite multiplicity, where $P'_i$ is some projection with $\pi(P'_i)=P_i$. we have 
relation (5).

Now we consider the universal $C^\ast$-algebra
generated by the family $\{ S_g , U_h \mid g\in\bigcup_{i\in
I}G_i\setminus H, h\in H\}$ satisfying (1), (2), (3) and (4). We denote it by $\mathcal{O}_{\Gamma}$. Here, the universality means that if another family 
$\{s_g, u_h\}$ satisfies (1), (2), (3) and (4), then there
exists a surjective $\ast$-homomorphism $\phi$ of
$\mathcal{O}_{\Gamma}$ onto $C^\ast(s_g, u_h)$ such that
$\phi(S_g)=s_g$ and $\phi(U_h)=u_h$. Summing up the above, we employ
the following definitions and notation:

\begin{definition}
Let $I$ be a finite index set and
$G_i$ be a group containing a copy of a finite group $H$ as a
subgroup for $i\in I$. Suppose that each $G_i$
is either a finite group or $\mathbb{Z}\times H$. Let $I_0$ be the subset of $I$ such
that $G_i$ is finite for all $i\in I_0$. We denote the
amalgamated free product $\ast_H G_i$ by $\Gamma$.

We fix a set $\Omega_i$ of left
representatives of $G_i/H$ with $e\in\Omega_i$ and a set $X_i$ of 
representatives of $H\backslash G_i/H$ which is
contained in $\Omega_i$. Let $(a_i, e)$ be a generator
of $G_i$ for $i\in I\setminus I_0$. We write
$a_i$, for short. Here we choose
$\Omega_i=X_i=\{a_i^n \mid n\in\mathbb{N} \}$.
We exclude the case where
$\bigcup_i\Omega_i\setminus\{e\}$ has only one or two points.

We define the corresponding universal $C^\ast$-algebra $\mathcal{O}_{\Gamma}$
generated by partial isometries $S_g$ for $g\in\bigcup_{i\in I}G_i\setminus H$ and unitaries $U_h$ for $h\in H$ satisfying (1), (2), (3) and (4).

We set for $\gamma\in\Gamma$,
$$\begin{array}{cc}
Q_{\gamma}=S_{\gamma}^{\ast}\cdot S_{\gamma}, & P_{\gamma}=S_{\gamma}\cdot 
S_{\gamma}^{\ast},
\end{array}$$
$$\begin{array}{cc}
P_i=\sum_{g\in\Omega_i}P_g & \mbox{if}\ i\in I_0, \\
P_i=P_{a_i}+P_{a_i^{-1}} & \mbox{if}\ i\in
I\setminus I_0.
\end{array}$$
For convenience, we set for any integer $n$,
$$\Gamma_n=\{\gamma\in\Gamma \mid |\gamma|=n\},$$ 
$$\Delta_n=\{\gamma\in\Gamma_n \mid \gamma=\gamma_1\cdots\gamma_n,
\gamma_k\in\Omega_{i_k}, i_1\ne\dots\ne i_n\}.$$
We also set $\Delta=\bigcup_{n\geq 1}\Delta_n$.
\end{definition}

\begin{lemma}
For $i\in I$ and $h\in H$,
$$U_h P_i=P_i U_h.$$
\end{lemma}
{\it Proof.}
Use the above relations (2). \hfill $\Box$

\begin{lemma}
Let $\gamma_1,\gamma_2 \in\Gamma$. Suppose that
$S_{\gamma_1}^{\ast}S_{\gamma_2}\ne 0.$

If $|\gamma_1|=|\gamma_2|$, then $S_{\gamma_1}^{\ast}S_{\gamma_2}=Q_g U_h \,$
 for some $g\in\bigcup_{i\in I}G_i, h\in H$.

If $|\gamma_1|>|\gamma_2|$,
then $S_{\gamma_1}^{\ast}S_{\gamma_2}=S_{\gamma}^{\ast}$ for some
$\gamma\in\Gamma$ with $|\gamma|=|\gamma_1|-|\gamma_2|$.

If $|\gamma_1|<|\gamma_2|$, then $S_{\gamma_1}^{\ast}S_{\gamma_2}=S_{\gamma}$
 for some $\gamma\in\Gamma$ with $|\gamma|=|\gamma_2|-|\gamma_1|$.
\end{lemma}
{\it Proof.}
By (2), we obtain the lemma. \hfill $\Box$

\begin{corollary}
$$\mathcal{O}_{\Gamma}=\overline{\rm span}\{ \, S_{\mu} P_i S_{\nu}^{\ast}
\mid \mu,\nu\in\Gamma, i\in I \, \}.$$
\end{corollary}
{\it Proof.}
This follows from the previous lemma. \hfill $\Box$

\medskip 

Next we consider the gauge action of $\mathcal{O}_{\Gamma}$. Namely, if
$z\in\mathbb{T}$ then the family $\{ \, zS_g, U_h \,\}$ also satisfies (1), (2), (3), (4) and generates $\mathcal{O}_{\Gamma}$. The
universality gives an automorphism $\alpha_z$ on $\mathcal{O}_{\Gamma}$
 such that $\alpha_z(S_g)=zS_g$ and $\alpha_z(U_h)=U_h$. In fact,
 $\alpha$ is a continuous action of 
$\mathbb{T}$ on
$\mathcal{O}_{\Gamma}$, which is called {\it the gauge action}. Let $dz$ be the
normalized Haar measure on $\mathbb{T}$ and we define a conditional
expectation $\Phi$ of $\mathcal{O}_{\Gamma}$ onto the fixed-point
algebra $\mathcal{O}_{\Gamma}^{\mathbb{T}}=\{ \,
a\in\mathcal{O}_{\Gamma} \mid \alpha_z(a)=a, \mbox{for}\ z\in\mathbb{T} \, \}$ by

$$\Phi(a)=\int_{\mathbb{T}}\alpha_z(a)\,dz, \qquad \mbox{for}\ \,
a\in\mathcal{O}_{\Gamma}.$$

\begin{lemma}
The fixed-point algebra $\mathcal{O}_{\Gamma}^{\mathbb{T}}$ is an AF-algebra.
\end{lemma}
{\it Proof.}
For each $i\in I$, set
$$ \mathcal{F}_n^i=\overline{\rm span}\{ \, S_{\mu}P_iS_{\nu}^\ast \mid
\mu, \nu\in\Gamma_n \, \}. $$
We can find systems of matrix units in
$\mathcal{F}_n^i$, parameterized by  $\mu, \nu\in\Delta_n$, as follows:
$$ e_{\mu, \nu}^i=S_{\mu}P_iS_{\nu}^\ast.$$
Indeed, using the previous lemma, we compute
$$ e_{\mu_1, \nu_1}^i e_{\mu_2, \nu_2}^i=\delta_{\nu_1,
  \mu_2}S_{\mu_1}P_i Q_{\nu_1}P_i S_{\nu_2}^\ast=\delta_{\nu_1,
  \mu_2} e_{\mu_1, \nu_2}^i.$$
Thus we obtain the identifications
$$ \mathcal{F}_n^i\simeq M_{N(n, i)}(\mathbb{C})\otimes e_{\mu, 
\mu}^i\mathcal{F}_n^i e_{\mu, \mu}^i, $$
for some integer $N(n, i)$ and some $\mu\in\Delta_n$. Moreover, for
$\xi, \eta$,
$$e_{\mu, \mu}^i \left(S_{\xi} P_iS_{\eta}^\ast \right) e_{\mu,
  \mu}^i = \left\{\begin{array}{ll}
    S_{\mu}P_i U_h P_i S_{\mu}^\ast & \mbox{if} \quad \xi, \eta\in\mu H, \\
                                 0 & \mbox{otherwise}.\
                          \end{array}\right.$$ 
for some $h\in H$. Note that $C^\ast(S_{\mu}P_iU_hP_iS_{\mu}^\ast\mid h\in H)$ 
is isomorphic to $C^\ast(P_iU_hP_i\mid h\in H)$ via the map $x\mapsto 
S_{\mu}^\ast xS_{\mu}$. Therefore the relation (5) gives
$$ \mathcal{F}_n^i\simeq M_k(\mathbb{C})\otimes\overline{\rm span}\{ \,
S_{\mu}P_iU_h P_i S_{\mu}^\ast \mid h\in H \, \}\simeq
M_k(\mathbb{C})\otimes C^\ast(H).$$
Note that $\{ \mathcal{F}_n^i \,\mid i\in I \, \}$ are mutually orthogonal and
$$ \mathcal{F}_n=\oplus_{i\in I}\mathcal{F}_n^i$$
is a finite-dimensional $C^\ast$-algebra.

The relation (2) gives
$\mathcal{F}_n\hookrightarrow\mathcal{F}_{n+1}$. Hence,
$$ \mathcal{F}=\overline{\bigcup_{n\geq 0}\mathcal{F}_n} $$
is an $AF$-algebra. Therefore it suffices to show that
$\mathcal{F}=\mathcal{O}_{\Gamma}^{\mathbb{T}}$. It is trivial that
$\mathcal{F}\subseteq\mathcal{O}_{\Gamma}^{\mathbb{T}}$. On the other
hand, we can approximate any $a\in\mathcal{O}_{\Gamma}^{\mathbb{T}}$
by a linear combination of elements of the form $S_{\mu}P_iS_{\nu}^\ast$. Since
$\Phi(a)=a$, $a$ can be approximated by a linear combination
of elements of the form $S_{\mu}P_iS_{\nu}^\ast$ with $|\mu|=|\nu|$. Thus
$a\in\mathcal{F}$. \hfill $\Box$

\medskip

We need another lemma to prove the uniqueness of $\mathcal{O}_{\Gamma}$.
\begin{lemma}
Suppose that $i_0\in I$ and $W$ consists of finitely many elements $(\mu, h)\in\Delta\times H$ 
such that the last word of $\mu$ is not contained in $\Omega_{i_0}$ and $W\cap
H=\emptyset$. Then there exists $\gamma=g_0\cdots g_n$ with
$g_k\in\Omega_{i_k}$ and $i_0\ne\dots\ne i_n\ne i_0$ such that for any $(\mu,
h)\in W$, $\mu h\gamma$ never have the form $\gamma\gamma'$ for some
$\gamma'\in\Gamma$. \end{lemma} {\it Proof.} Let $i_0\in I$ and $W$ be a
finite subset of $\Delta\times H$ as above. We first assume that $|I|\geq 3$.
Then we can choose $x\in\Omega_{i_0}, y\in\Omega_{j}$ and $z\in\Omega_{j'}$
such that $j\ne i_0\ne j'$ and $j\ne j'$. For sufficiently long word 
$$\gamma=(xy)(xz)(xyxy)(xzxz)(xyxyxy)(xzxzxz)\cdots(\cdots z),$$
we are done. We next assume that $|I|=2$. Since we exclude the case where
$\Omega_1\cup\Omega_2\setminus\{e\}$ has only one or two elements, we
can choose at least three distinct points $x\in\Omega_{i_0},
y\in\Omega_{j}$ and $z\in\Omega_{j'}$. If $i_0\ne j=j'$ we set 
$$\gamma=(xy)(xz)(xyxy)(xzxz)(xyxyxy)(xzxzxz)\cdots(\cdots z),$$
as well. If $i_0=j\ne j'$ we set
$$\gamma=(xz)(yz)(xzxz)(yzyz)(xzxzxz)(yzyzyz)\cdots(\cdots z).$$
Then if $\gamma$ has the desired properties, we are done. Now assume that there exist
some $(\mu, h)\in W$ such that $\mu h\gamma=\gamma\gamma'$ for some
$\gamma'$. Fix such an element $(\mu, h)\in W$. By hypothesis, we can choose
$\delta\in\Delta$ with $|\gamma'|\leq|\delta|$ such that the last
word of $\delta$ does not belong to $\Omega_{i_0}$ and $\delta$ does not have the form
$\gamma'\delta'$ for some $\delta'$. 
Set
$\tilde{\gamma}=\gamma\delta$. Then $\mu h\tilde{\gamma}$ does not have the form
$\gamma\gamma''$ for any $\gamma''$. Indeed,
$$\mu h\tilde{\gamma}=\mu h\gamma\delta=\gamma\gamma'\delta\ne\tilde{\gamma}\gamma'',$$
for some $\gamma''$.
Since $W$ is finite, we can obtain a desired element $\gamma$ by replacing $\tilde{\gamma}$, inductively.
\hfill $\Box$
\medskip

We now obtain the uniqueness theorem for $\mathcal{O}_{\Gamma}$.

\begin{theorem}
Let $\{ \, s_g, u_h \,\}$ be another family of partial isometries and
unitaries satisfying (1), (2), (3) and (4).
Assume that 
$$C^\ast(H)\simeq C^\ast( \, p_iu_hp_i \mid h\in H \, ),$$
where $p_i=\sum_{g\in\Omega_i\setminus\{e\}}s_gs_g^\ast$ for
$i\in I_0$ and
$p_i=s_{a_i}s_{a_i}^\ast+s_{a_i^{-1}}s_{a_i^{-1}}^\ast
$ 
for $i\in I\setminus I_0$. Then the canonical surjective $\ast$-homomorphism $\pi$ of
$\mathcal{O}_{\Gamma}$ onto $C^\ast\left( \, s_g, u_h \, \right)$ is
faithful.
\end{theorem}
{\it Proof.}
To prove the theorem, it is enough to show that (a) $\pi$ is faithful on the fixed-point algebra
$\mathcal{O}_{\Gamma}^{\mathbb{T}}$, and (b)
$\|\pi\left(\Phi(a)\right)\|\leq\|\pi(a)\|$ for all
$a\in\mathcal{O}_{\Gamma}$ thanks to ~\cite[Lemma 2.2]{bkr}.

To establish (a), it suffices to show that $\pi$ is faithful on
 $\mathcal{F}_n$ for all $n\geq 0$. By the proof of Lemma 4.5, we have
$$ \mathcal{F}_n^i=M_{N(n, i)}(\mathbb{C})\otimes C^\ast(H), $$
for some integer $N(n, i)$. Note that $s_g s_g^\ast$ is 
non-zero. Hence $\pi$ is
injective on $M_{N(n, i)}(\mathbb{C})$. By the other hypothesis, $\pi$ is
injective on $C^\ast(H)$. 

Next we will show (b). It is enough to check (b) for
$$a=\sum_{\mu, \nu\in F}\sum_{j\in J}C_{\mu, \nu}^j S_{\mu} P_j 
S_{\nu}^\ast,$$
where $F$ is a finite subset of $\Gamma$ and $J$ is a subset of $I$. 
For $n=\max\{ |\mu| \mid \mu\in F \}$, we have
$$ \Phi(a)=\sum_{\{ \mu, \nu\in F \mid |\mu|=|\nu| \}}\sum_{j\in J}C_{\mu, 
\nu}^j
S_{\mu} P_j S_{\nu}^\ast\in\mathcal{F}_n. $$

Now by changing $F$ if necessary, we may assume that $\min\{ |\mu|, 
|\nu| \}=n$ for every pair $\mu, \nu\in F$ with $C_{\mu, \nu}^j\ne
0$. Since $\mathcal{F}_n=\oplus_i\mathcal{F}_n^i$,
there exists some $i_0\in J$ such that 
$$ \|\pi(\Phi(a))\|=\|\sum_{|\mu|=|\nu|}C_{\mu, \nu}^{i_0}
s_{\mu}p_{\imath_0}s_{\nu}^\ast\|.$$

By changing $F$ such that $F\subset\Delta$ again, we may further assume that
$$ \|\pi(\Phi(a))\|=\|\sum_{\stackrel{\mbox{$\scriptstyle \mu, \nu\in
      F$}}
                                     {\mbox{$\scriptstyle
      |\mu|=|\nu|$}}}
                      \sum_{h\in F'}C_{\mu, \nu, h}^{i_0}
      s_{\mu}p_{i_0}u_hp_{i_0}s_{\nu}^\ast\|$$
where $F'$ consists of elements of $H$, (perhaps with multiplicity). By applying the preceding lemma to 
$$W=\{(\mu', h)\in\Delta\times H \mid \mu'\enskip\mbox{is subword of}\enskip\mu\in F, h^{-1}\in F'\},$$
we have $\gamma\in\Delta$ satisfying the property in the previous
lemma. Then 
we define a projection
$$
Q=\sum_{\tau\in\Delta_n}s_{\tau}s_{\gamma}p_{i_0}s_{\gamma}^\ast 
s_{\tau}^\ast.$$
By hypothesis, $Q$ is non-zero.

If $\mu, \nu\in\Delta_n$ then
$$
Q\left(s_{\mu}p_{i_0}s_{\nu}^\ast\right)Q=s_{\mu}s_{\gamma}p_{i_0}s_{\gamma}^\ast 
p_{i_0}s_{\gamma}p_{i_0}s_{\gamma}^\ast s_{\nu}^\ast=s_{\mu}s_{\gamma}p_{i_0}s_{\gamma}
^\ast 
s_{\nu}^\ast$$
is non-zero. Therefore $s_{\mu}(s_{\gamma}p_{i_0}s_{\gamma}^\ast)s_{\nu}^\ast$ is 
also a family of matrix units
  parameterized by $\mu, \nu\in\Delta_n$. Hence the same arguments as in the 
  proof of Lemma 4.5 give
$$ \pi(\mathcal{F}_n^{i_0})\simeq M_{N(n, i_0)}(\mathbb{C})\otimes C^\ast\left
  ( \, s_{\mu}s_{\gamma}p_{i_0}u_hp_{i_0}s_{\gamma}^\ast
  s_{\mu}^\ast \mid h\in H \, \right). $$
By hypothesis, we deduce that $b\mapsto Q\pi(b)Q$ is faithful on
$\mathcal{F}_n^{i_0}$. In particular, we conclude that
$\|\pi(\Phi(a))\|=\|Q\pi(\Phi(a))Q\|$.

We next claim that $Q\pi(\Phi(a))Q=Q\pi(a)Q$. We fix $\mu, \nu\in F$. If $|\mu|\ne|\nu|$ 
then one of $\mu, \nu$ has length $n$ and
the other is longer; say $|\mu|=n$ and $|\nu|>n$. Then
$$
Q\left(s_{\mu}p_{i_0}u_hp_{i_0}s_{\nu}^\ast\right)Q=s_{\mu}s_{\gamma}p_{i_0}s_{\gamma}^\ast 
p_{i_0}u_hp_{i_0}s_{\nu}^\ast\left(\sum_{\tau\in\Delta_n}s_{\tau}s_{\gamma}p_{i_0}s_{\gamma}^\ast 
  s_{\tau}^\ast\right).$$
Since $|\nu|>|\tau|$, this can have a non-zero summand only if
$\nu=\tau{\nu}'$ for some $\nu'$. However $ s_{\gamma}^\ast u_hs_{\nu}^\ast
s_{\tau}s_{\gamma}=s_{\gamma}^\ast u_hs_{{\nu}'}^\ast s_{\gamma}$, and $s_{{\nu}'h^{-1}\gamma}^\ast s_{\gamma}$ is non-zero only if
${\nu}'h^{-1}\gamma$ has the form $\gamma\gamma'$. This is impossible by the choice of
$\gamma$. Therefore we have $Q\left(s_{\mu}p_{i_0}s_{\nu}\right)Q=0$ if
  $|\mu|\ne|\nu|$, namely $Q\pi(\Phi(a))Q=Q\pi(a)Q$. Hence we can finish 
  proving (b):
$$\|\pi(\Phi(a))\|=\|Q\pi(\Phi(a))Q\|=\|Q\pi(a)Q\|\leq\|\pi(a)\|.$$
Therefore ~\cite[Lemma 2.2]{bkr} gives the theorem. \hfill $\Box$

\medskip

By essentially the same arguments, we can prove the following.

\begin{corollary}
Let $\{t_g, v_h\}$ and $\{s_g, u_h\}$ be two families of partial isometries and unitaries satisfying (1), (2), (3) and (4). Suppose that the map $p_iv_hp_i\mapsto q_iu_hq_i$ gives an isomorphism:
$$C^\ast(p_iv_hp_i\mid h\in H)\simeq C^\ast(q_iv_hq_i\mid h\in H),$$
where $p_i=\sum_{g\in \Omega_i\setminus\{e\}}t_gt_g^\ast, q_i=\sum_{g\in \Omega_i\setminus\{e\}}s_gs_g^\ast$ and so on. Then the canonical map gives the isomorphism between $C^\ast(t_g, v_h)$ and $C^\ast(s_g, u_h)$. 
\end{corollary}

\medskip

Before closing this section, we will show that our algebra $\mathcal{O}_{\Gamma}$ is
isomorphic to a
certain Cuntz-Krieger-Pimsner algebra. Let $A=C^\ast\left( \, P_iU_hP_i \mid 
h\in H, i\in I \, 
\right)\simeq \bigoplus_{i\in I}C_r^\ast(H)$. We define a Hilbert $A$-bimodule $X$ 
as follows:
$$ X=\overline{\rm span}\{ \, S_gP_{i} \mid g\in\bigcup_{j\ne i}G_j, \, 
|g|=1, \, i\in I \,
  \}$$
with respect to the inner product $\langle S_gP_i,
  S_{g'}P_j\rangle=P_iS_g^\ast S_{g'}P_j \in A.$
In terms of the groups, the $A$-$A$ bimodule structure can be described as follows: we set 
$$A=\bigoplus_{i\in I}A_i=\bigoplus_{i\in I}\mathbb{C}[H],$$
and define an $A$-bimodule $\mathcal{H}_i$ by
$$\mathcal{H}_i=\mathbb{C}[\{g\in\bigcup_{j\ne i}G_j \mid |g|=1\}]$$ 
with left and right $A$-multiplications such that for $a=(h_i)_{i\in
  I}\in A$ and $g\in G_j\setminus H\subset \mathcal{H}_i$,
$$a\cdot g=h_jg \quad\mbox{and}\quad g\cdot a=gh_i,$$
and with respect to the inner product 
$$\langle g, g' \rangle_{\mathcal{H}_i}=\left\{\begin{array}{ll}
   g^{-1}g'\in A_i   & \mbox{if}\ g^{-1}g'\in H, \\ 
   0   & \mbox{otherwise}.
   \end{array}\right. $$
Then we define the $A$-bimodule $X$ by 
$$X=\bigoplus_{i\in I}\mathcal{H}_i,$$
and we obtain the CKP-algebra $\mathcal{O}_X$.

\begin{proposition} 
Assume that $A$ and $X$ are as above. Then
$$\mathcal{O}_{\Gamma}\simeq\mathcal{O}_X.$$
\end{proposition}
{\it Proof.}
We fix a finite basis $u(g, i)=g\in\mathcal{H}_i$ for
  $g\in\Omega_j, i\in I$ with $j\ne i, |g|=1$. Then we have
  $\mathcal{O}_X=C^\ast(S_{u(g, i)})$. Let $s_{u(g, i)}=S_gP_i$ in $\mathcal{O}_{\Gamma}$. Note
  that we have $\mathcal{O}_{\Gamma}=C^\ast(s_{u(g, i)})$. The
  relation (4) corresponds to the
  relations $(\dag)$ of the CKP-algebras. The
  family $\{s_{u(g, i)}\}$ therefore satisfies the relations of
  the CKP-algebras. Since the CKP-algebra has universal properties, there
  exists a canonical surjective $\ast$-homomorphism of $\mathcal{O}_X$ onto
$\mathcal{O}_{\Gamma}$. Conversely, let $s_g=\sum_{i\in I}S_{u(g, i)}$ and
$u_h=\oplus_{i\in I}h$ for $h\in H$ in $\mathcal{O}_X$, and then we have $\mathcal{O}_X=C^\ast(s_g, u_h)$. 
By the universality of
$\mathcal{O}_{\Gamma}$, we can also obtain a canonical surjective
$\ast$-homomorphism of $\mathcal{O}_{\Gamma}$ onto
  $\mathcal{O}_X$. These maps are mutual inverses. Indeed,
$$\left. \begin{array}{ccccc}
   S_g & \mapsto & \sum_{i\in I}S_{u(g, i)} & \mapsto &
   \sum_{i\in I}S_gP_i=S_g, \\
   U_h & \mapsto & \bigoplus_{i\in I}h & \mapsto &
   \sum_{i\in I}P_iU_hP_i=U_h.
\end{array}\right. $$ \hfill $\Box$

\section{Crossed product algebras associated with
$\mathcal{O}_{\Gamma}$}

In this section, we will show that $\mathcal{O}_{\Gamma}$ is isomorphic to a crossed product
algebra. We first define a ``boundary space''. We set
$$\tilde{\Lambda}=\{ \, (\gamma_n)_{n\geq 0} \mid \gamma_n \in\Gamma,
|\gamma_n|+1=|\gamma_{n+1}|, |\gamma_n^{-1}\gamma_{n+1}|=1 \,
\mbox{for a sufficiently large}\ n\geq 0 \, \}.$$
We introduce the following equivalence relation $\sim$; $(\gamma_n)_{n\geq
  0}, (\gamma'_n)_{n\geq 0}\in\tilde{\Lambda}$ are equivalent if there exists some 
$k\in\mathbb{Z}$ such that $\gamma_n H =
\gamma'_{n+k} H$ for a sufficiently large $n$. Then we define
$\Lambda=\tilde{\Lambda}/\sim$. We denote the equivalent class of
$(\gamma_n)_{n\geq 0}$ by $[\gamma_n]_{n\geq 0}$. 

Before we define an action of $\Gamma$ on $\Lambda$, we construct
another space $\Omega$ to introduce a compact space structure, on which $\Gamma$ acts continuously. Let $\Omega$ denote the set of
sequences $x : \mathbb{N}\to\Gamma$ such that 
$$ \left\{\begin{array}{ll}
   x(n)\in\Omega_{i_n}\setminus\{e\} & \mbox{for}\ n\geq 1,
   \\
   x(n)\in\{a_{i_n}^{\pm 1}\} & \mbox{if}\ i_n \in
   I\setminus I_0, \\
   i_n \ne i_{n+1} & \mbox{if}\ i_n \in I_0, \\
   x(n)=x(n+1) & \mbox{if}\ i_n \in I\setminus I_0,
   i_n=i_{n+1}. 
\end{array}\right. $$
Note that $\Omega$ is a compact Hausdorff subspace of 
$\prod_{\mathbb{N}}\left(\bigcup_i\Omega_i\setminus\{e\}\right)$. We introduce a map $\phi$ between $\Lambda$ and
$\Omega$; for $x=(x(n))_{n\geq 1}\in\Omega$, we define a map
$\phi(x)=[\gamma_n]\in\Lambda$ by
\begin{eqnarray*}
\gamma_0 &=& e \qquad \mbox{if}\ n=0, \\
\gamma_n &=& x(1)\cdots x(n), \qquad \mbox{if}\ n\geq 1.
\end{eqnarray*}

\begin{lemma}
The above map $\phi$ is a bijection from $\Lambda$ onto $\Omega$ and hence $\Lambda$ inherits a compact space structure via $\phi$. 
\end{lemma}
{\it Proof.}
For
$x=(x(n))\ne{x}'=({x}'(n))$, there exists an
integer $k$ such that $x(k)\ne{x}'(k)$. If
$\phi(x)=[\gamma_n]$ and $\phi({x}')=[{\gamma}'_n]$, then
$\gamma_k H\ne{\gamma}'_k H$. Hence we have injectivity of
$\phi$. Next we will show surjectivity. Let
$[\gamma_n]\in\Sigma$. We may take a representative $(\gamma_n)$ satisfying $|\gamma_n|=n$. 
Now we assume that $\gamma_n$ is uniquely expressed as $\gamma_n=g_1\cdots g_n h, \gamma_{n+1}
=g'_1\cdots g'_{n+1}h'$ for
$g_k\in\Omega_{i_k}, g'_k\in\Omega_{j_k}, h, h'\in
H$. Since $|\gamma_n^{-1}\gamma_{n+1}|=1$, we have
$$ h^{-1}g_n^{-1}\cdots g_1^{-1}g'_1\cdots g'_{n+1}h'=g, $$
for some $g\not\in H$ with $|g|=1$. Inductively, we have $g_1=g'_1, \dots ,
g_n=g'_n$. Hence we can assume that
$\gamma_n=g_1\cdots g_n$. We set $x(n)=g_n$ and get
$\phi((x(n)))=[\gamma_n]$. \hfill $\Box$

\medskip

Next we define an action of $\Gamma$ on $\Lambda$. Let
$[\gamma_n]_{n\geq 0}\in\Lambda$. For $\gamma\in\Gamma$, define 
$$\gamma\cdot[\gamma_n]_{n\geq 0}=[\gamma\gamma_n]_{n\geq 0}.$$
We will show that this is a continuous action of $\Gamma$ on $\Lambda$. Let $[\gamma_n]$,$[\gamma'_n]\in\Lambda$ such that 
$(\gamma_n)\sim
(\gamma'_n)$ and $\gamma\in\Gamma$. Since there exists some integer $k$ such that $\gamma_n
H=\gamma'_{n+k} H$ for sufficiently large integers $n$, we have $\gamma\gamma_n
H=\gamma\gamma'_{n+k} H$. Hence this is
well-defined. To show that $\gamma$ is continuous, we consider how $\gamma$ acts on $\Omega$ 
via the map $\phi$. For $g\in \Omega_{i}$ with $|g|=1$ and $x=(x(n))_{n\geq 1}\in\Omega$,
$$ (g\cdot x)(1)=\left\{\begin{array}{ll}
   g   & \mbox{if}\ i\ne i_1, \\   
   g_1 & \mbox{if}\ i=i_1, \, gx(1)\not\in H, \,
   i\in I_0, \\
       & \mbox{and}\ gx(1)=g_1h_1 \, (g_1\in\Omega_{i_1},
   h_1\in H ), \\
   g   & \mbox{if}\ i=i_1, \, gx(1)\not\in H, \,
   i\in I\setminus I_0, \\
   g_2 & \mbox{if}\ i=i_1, \, gx(1)\in H, \, i\in 
   I_0, \\
       & \mbox{and}\ gx(1)=h_1, \, h_1x(2)=g_2h_2
    (g_2\in\Omega_{i_2}, h_1, h_2\in H ), \\
   x(2)    & \mbox{if}\ i=i_1, \, gx(1)\in H, \,
   i\in I\setminus I_0, 
   \end{array}\right. $$
and for $n>1$,
$$ (g\cdot x)(n)=\left\{\begin{array}{ll}
   x(n-1)   & \mbox{if}\ i\ne i_1, \\
   g_n            & \mbox{if}\ i=i_1, \, gx(1)\not\in H, \\
                  & \mbox{and}\ h_{n-1}x(n)=g_nh_n \, (g_n\in\Omega_{i_n},
   h_n\in H ), \\
   x(n-1)   & \mbox{if}\ i=i_1, \, gx(1)\not\in H, \,
   i\in I\setminus I_0, \\
   g_{n+1}        & \mbox{if}\ i=i_1, \, gx(1)\in H, \\
                  & \mbox{and}\ h_nx(n+1)=g_{n+1}h_{n+1},
   (g_{n+1}\in\Omega_{i_{n+1}}, h_{n+1}\in H ), \\
   x(n+1)   & \mbox{if}\ i=i_1, \, gx(1)\in H, \,
   i\in I\setminus I_0.
   \end{array}\right. $$

For $h\in H$,
$$ (h\cdot x)(n)=\left\{\begin{array}{ll}
   g_1   & \mbox{if}\ n=1, \\
         & \mbox{and}\ hx(1)=g_1h_1, \,
   (g_1\in\Omega_{i_1}, h_n\in H), \\  
   g_n   & \mbox{if}\ n>1, \\
         & \mbox{and}\ h_{n-1}x(n)=g_nh_n, \,
   (g_n\in\Omega_{i_n}, h_n\in H).
   \end{array}\right. $$

Then one can check easily that the pull-back of any open set of $\Omega$ by $\gamma$ is also an open set of $\Omega$. Thus we have proved that $\gamma$ is a homeomorphism on $\Lambda$. The 
equations
$$(\gamma\gamma')[\gamma_n]=[\gamma\gamma'\gamma_n]=\gamma([\gamma'\gamma_n])=\gamma\circ\gamma'[\gamma_n],$$
imply associativity.

Therefore we have obtained the following:
\begin{lemma}
The above space $\Omega$ is a compact Hausdorff space and $\Gamma$
acts on $\Omega$ continuously. 
\end{lemma}

The following result is the main theorem of this section.

\begin{theorem}
Assume that $\Omega$ and the action of $\Gamma$ on $\Omega$ are as above. Then we have the 
identifications
$$\mathcal{O}_{\Gamma}\simeq C(\Omega)\rtimes\Gamma\simeq C(\Omega)\rtimes_r\Gamma.$$
\end{theorem}
{\it Proof.}
We first consider the full crossed product $C(\Omega)\rtimes\Gamma$. Let
$Y_{i}=\{ \, (x(n)) \mid x(1)\in\Omega_{i} \, \}\subset\Omega$
 be clopen sets for ${i}\in I$. Note that if $i\in I_0$,
 then $Y_{i}$ is the
disjoint union of the clopen sets $\{ \, g(\Omega\setminus Y_{i})
\mid g\in\Omega_i\setminus\{e\} \, \}$, and if $i\in
I\setminus I_0$, then $Y_{i}=Y_{i}^{+}\cup Y_{i}^{-}$
where $Y_i^{\pm}=\{ \, (x(n)) \mid x(1)=a_i^{\pm} \, \}$. Let
$p_i=\chi_{\Omega\setminus Y_i}$ and
$p_i^{\pm}=\chi_{Y_i^{\pm}}$. We define $T_g=gp_i$ for $g\in
G_i\setminus H$ and $i\in I_0$ and $T_{a_i^{\pm
    1}}=a_i^{\pm 1}\left(p_i+p_i^{\pm}\right)$ for $i\in 
I\setminus I_0$. Let
$V_h=h$ for $h\in H$. Then the family $\{ T_g, V_h \}$ satisfies the
relations (1), (2), (3) and (4). Indeed, we can first check that $h\in H$ commutes with
$p_i$ and $p_i^{\pm 1}$. So the relation
(1) holds. Let $g\in G_i\setminus H$ and $g'\in G_j\setminus H$ 
with $i, j\in I_0$. Then 
$$T_g^\ast
T_{g'}=p_ig^{-1}g'p_j=g^{-1}\chi_{g(\Omega\setminus
  Y_i)}\chi_{g'(\Omega\setminus Y_j)}g'=\delta_{i, j}\delta_{gH, g'H}p_ig^{-1}g'.$$
Moreover it follows from $\Omega\setminus
Y_i=\bigcup_{j\ne i}Y_j$ that
\begin{eqnarray*}  
\lefteqn{T_g^\ast T_g=\chi_{\Omega\setminus 
Y_i}=\sum_{j\ne i}\chi_{Y_j}} \\
 &=& \sum_{j\in
 I_0, j\ne i}\sum_{g\in\Omega_j\setminus\{e\}}\chi_{g(\Omega\setminus
 Y_j)}+\sum_{j\in I\setminus I_0}\chi_{a_j(\Omega\setminus
 Y_j)}+\chi_{a_j^{-1}(\Omega\setminus Y_j)} \\
 &=& \sum_{j\in
 I_0, j\ne i}\sum_{g\in\Omega_j\setminus\{e\}}gp_jg^{-1}+\sum_{j\in I\setminus 
I_0}p_j^++p_j^{-} \\
 &=& \sum_{j\in
 I_0, j\ne i}\sum_{g\in\Omega_j\setminus\{e\}}T_gT_g^\ast+\sum_{j\in I\setminus 
I_0}T_{a_j}T_{a_j}^\ast+T_{a_j^{-1}}T_{a_j^{-1}}^\ast.
\end{eqnarray*}
For all other cases, we can also check the relations (2) and (3) by similar calculations. Since $\Omega$ is the disjoint union of $Y_i$, we have (4). Note that $g, 
p_i, p_i^{\pm}\in C^\ast(T_g, V_h)$. Moreover, since the family 
$\{\gamma(\Omega\setminus Y_i)\mid \gamma\in\Gamma, i\in
I\}\cup\{\gamma Y_i^{\pm}\mid \gamma\in\Gamma, i\in
I\setminus I_0\}$ generates the topology of $\Omega$, we have
$C(\Omega)\rtimes\Gamma=C^\ast(T_g, V_h)$. By the
universality of $\mathcal{O}_{\Gamma}$, there exists a canonical
surjective $\ast$-homomorphism of $\mathcal{O}_{\Gamma}$ onto
$C(\Omega)\rtimes\Gamma$, sending $S_g$ to $T_g$ and $U_h$ to $V_h$. 

Conversely, let $q_i=\sum_{j\ne i}P_j$ and
$q_i^{\pm}=S_{a_i^{\pm 1}}S_{a_i^{\pm 1}}^\ast$. Let
$$ \left\{\begin{array}{ll}
w_g=S_g+\sum_{g'\in\Omega_i\setminus H\cup g^{-1}H}S_{gg'}S_{g'}^\ast+S_{g}^{\ast}
& \mbox{for}\ g\in G_i\setminus H, i\in I_0, \\
w_{a_i}=S_{a_i}+S_{a_i^{-1}}^\ast & \mbox{for}\
i\in I\setminus I_0, \\
w_h=U_h & \mbox{for}\ h\in H.
\end{array}\right. $$
We will check that $w_g$ are unitaries for $g\in
G_i\setminus H$ with $i\in I_0$. If
$g'\in\Omega_i\setminus H\cup g^{-1}H$, then
$gg'H=\gamma H$ for some $\gamma\in\Omega_i\setminus\{e, g\}$. Hence 
\begin{eqnarray*}
\lefteqn{w_gw_g^\ast} \\
 &=& \left(S_g+\sum_{g'\in\Omega_i\setminus H\cup
 g^{-1}H}S_{gg'}S_{g'}^\ast+S_{g^{-1}}^{\ast}\right)\left(S_g+\sum_{g'\in\Omega_i\setminus H\cup g^{-1}H}S_{gg'}S_{g'}^\ast+S_{g^{-1}}^{\ast}\right)^\ast \\
 &=& S_gS_g^\ast+\sum_{g'\in\Omega_i\setminus H\cup g^{-1}H}S_{gg'}S_{g'}^\ast S_{g'}S_{gg'}^\ast+S_{g^{-1}}^\ast
      S_{g^{-1}} \\
 &=& P_g+\sum_{g'\in\Omega_i\setminus\{e, g\}}P_{g'}+Q_g=1.
\end{eqnarray*}
Similarly, we have $w_g^\ast w_g=1$. For the other case, we can check
in the same way.

If $i\in I_0, \tau\in\Omega_i\setminus\{e\}$ then
\begin{eqnarray*}
\lefteqn{\sum_{g\in\Omega_i}w_gq_iw_g^\ast} \\
 &=&
 \sum_{g\in\Omega_i}\left(S_g+\sum_{g'\in\Omega_i\setminus H\cup g^{-1}H}S_{gg'}S_{g'}^\ast+S_{g^{-1}}^\ast\right)S_{\tau}^\ast S_{\tau}w_g^\ast \\
 &=&
 \sum_{g\in\Omega_i}S_gS_{\tau}^\ast 
S_{\tau}\left(S_g^\ast+\sum_{g'\in\Omega_i\setminus H\cup g^{-1}H}S_{g}S_{gg'}^\ast+S_{g^{-1}}\right) \\
 &=& \sum_{g\in\Omega_i}S_gS_{\tau}^\ast S_{\tau}S_g^\ast=1.
\end{eqnarray*}
For $i\in I\setminus I_0$, we have
$q_i^{+}+w_{a_i}q_i^{-}w_{a_i}^\ast=1$ and 
$q_i^{+}+q_i^{-}+q_i=1$ as well. Therefore the conjugates 
of the family $\{q_i, q_i^{\pm}\}$ by the elements of
$\Gamma$ generate a commutative $C^\ast$-algebra. This is the image of 
a representation of $C(\Omega)$. Therefore $(q_i, w)$ gives a covariant representation of the
$C^\ast$-dynamical system $(C(\Omega), \Gamma)$. Note that
$(q_i, w_g)$ generates $\mathcal{O}_{\Gamma}$. Hence by the universality of the full 
crossed product
$C(\Omega)\rtimes\Gamma$, there exists a canonical surjective
$\ast$-homomorphism of $C(\Omega)\rtimes\Gamma$ onto
$\mathcal{O}_{\Gamma}$. It is easy to show that the above two
$\ast$-homomorphisms are the
inverses of each other.

$$\begin{array}{ccccc}
S_g & \mapsto & gp_i & \mapsto & w_gQ_g=S_g, \\
S_{a_i^{\pm 1}} & \mapsto & a_i^{\pm
  1}(p_i+p_i^{\pm}) & \mapsto & w_{a_i^{\pm
  1}}(Q_{a_i^{\pm 1}}+P_{a_i^{\pm
  1}})=S_{a_i^{\pm 1}}, \\
U_h & \mapsto & h & \mapsto & U_h.
\end{array}$$

We have shown the identification $\mathcal{O}_{\Gamma}\simeq
C(\Omega)\rtimes\Gamma$. Since there exists a canonical surjective map of
$C(\Omega)\rtimes\Gamma$ onto $C(\Omega)\rtimes_r\Gamma$, we have a surjective $\ast$-homomorphism of
$\mathcal{O}_{\Gamma}$ onto $C(\Omega)\rtimes_r\Gamma$. Let 
$C(\Omega)\rtimes_r\Gamma=C^\ast(\tilde{\pi}(p_i), \lambda)$ where $\tilde{\pi}$ is 
the induced representation on the Hilbert space $l^2(\Gamma, \mathcal{H})$ by the universal 
representation $\pi$ of $C(\Omega)$ on a Hilbert space $\mathcal{H}$ and $\lambda$ is the 
unitary representation of $\Gamma$ on $l^2(\Gamma, \mathcal{H})$ such that 
$(\lambda_sx)(t)=x(s^{-1}t)$ for $x\in l^2(\Gamma, \mathcal{H})$. By the uniqueness theorem 
for $\mathcal{O}_{\Gamma}$, it
suffices to check 
$$
C^\ast\left(\tilde{\pi}(\chi_{Y_i})\lambda_h\tilde{\pi}(\chi_{Y_{i}})\right)\simeq C^\ast(H).$$
But the unitary representation
$\tilde{\pi}(\chi_{Y_i})\lambda_h\tilde{\pi}(\chi_{Y_{i}})$ is
quasi-equivalent to the left
regular representation of $H$. This completes the proof of the theorem. 
\hfill $\Box$

\medskip

In ~\cite{ser}, Serre defined the tree $G_T$, on which $\Gamma$
acts. In an appendix, we will give the definition of the tree $G_T=(V, E)$ where $V$ is the 
set of
vertices and $E$ is the set of edges. We denote the corresponding
natural boundary by $\partial G_T$. We also show how to construct
boundaries of trees in the appendix. (See Furstenberg 
~\cite{fur} and Freudenthal ~\cite{fre} for details.) 

\begin{proposition}
The space $\partial G_T$ is homeomorphic to $\Omega$ and the above two actions
of $\Gamma$ on $\partial G_T$ and $\Omega$ are conjugate.
\end{proposition}
{\it Proof.}
We define a map $\psi$ from $\partial G_T$ to
$\Omega$. First we assume that $I=\{1, 2\}$. The corresponding tree
$G_T$ consists of the vertex set $V=\Gamma/G_1\coprod\Gamma/G_2$
and the edge set $E=\Gamma/H$. For $\omega\in\partial G_T$, we can identify $\omega$ with an
infinite chain $\{G_{i_1}, g_1G_{i_2},
g_1g_2G_{i_3}, \ldots\}$ with
$g_k\in\Omega_{i_k}\setminus\{e\}$ and
$i_1\ne i_2\ne\cdots$. Then we define
$\psi(\omega)=[x(n)=g_{i_n}]$. We will recall the definition of the corresponding tree $G_T$, in general, on the appendix, (see ~\cite{ser}). Similarly, we can identify $\omega\in\partial G_T$ with an
infinite chain $\{G_0, G_{i_1}, g_1G_0,
g_1G_{i_2}, g_1g_2G_0,
\ldots\}$. Moreover we may ignore vertices $\gamma G_0$ for an infinite chain $\omega$, 
$$\{G_0, G_{i_1}, (g_1G_0\to\text{ignoring}),
g_1G_{i_2}, (g_1g_2G_0\to\text{ignoring}), g_1g_2G_{i_3}, \ldots\}.$$ Therefore, 
we define a map $\psi$ of $\partial G_T$ to $\Omega$ 
by
$$\psi(\omega)=[x(n)=g_{n}].$$
The pull-back by $\psi$ of any open set of $\partial G_T$ is an open set on
$\Omega$. It follows that $\psi$ is a homeomorphism. The two actions on $\partial G_T$ and 
$\Omega$ are defined by left
multiplication. So it immediately follows that these actions are conjugate.
\hfill $\Box$

\medskip

It is known that $\Gamma$ is a hyperbolic group (see a proof in
the appendix, where we recall the notion of
hyperbolicity for finitely generated groups as introduced by Gromov
 e.g. see ~\cite{gph}). Let $S=\{\bigcup_{i\in I}G_i\}$
and $G(\Gamma, S)$ be the Cayley graph of $\Gamma$ with the word
metric $d$. Let $\partial\Gamma$ be the hyperbolic boundary.

\begin{proposition}
The hyperbolic boundary $\partial\Gamma$ is homeomorphic to $\Omega$ and the actions of 
$\Gamma$ 
are conjugate.
\end{proposition}
{\it Proof.}
We can define a map $\psi$ from $\Omega$ to $\partial\Gamma$ by
$(x(n))\mapsto [x_n=x(1)\cdots x(n)]$. Indeed, since
$\langle x_n \, | \, x_m \rangle=\min\{n, m\}\to\infty \, (n,m\to\infty)$, it is
well-defined. For $x\ne y$ in $\Omega$, there exists $k$ such that $x(k)\ne
y(k)$. Then $\langle \psi(x) | \psi(y) \rangle\leq k+1$, which shows injectivity. Let
$(x_n)\in\partial\Gamma$. Suppose that $x_n=g_{n(1)}\cdots
g_{n(k_n)}h_n$ for some
$g_l\in\bigcup_i\Omega_i\setminus\{e\}$ with
$n(1)\ne\cdots\ne n(k_n)$. If $g_{n(1)}=g_{m(1)}, \ldots ,
g_{n(l)}=g_{m(l)}$ and $g_{n(l+1)}\ne g_{m(l+1)}$, then we set
$a_{n,m}=g_{n(1)}\cdots g_{n(l)}=g_{m(1)}\cdots g_{m(l)}$. So we have 
$$ \langle x_n \, | \, x_m \rangle\leq d(e,a_{n,m})+1\to\infty \, \, (n,m\to\infty). $$
Therefore we can choose sequences $n_1<n_2<\cdots$, and $m_1<m_2<\cdots$,
such that $a_{n_k, m_k}$ is a sub-word of $a_{n_{k+1}, m_{k+1}}$. Then
a sequence $\{ g_{n_k(1)}, \ldots , g_{n_k(l)}, g_{n_{k+1}(l+1)},
\ldots \}$ is mapped to $(x_n)$ by $\psi$. We have proved that $\psi$ is
surjective. The pull-back of any open set in $\partial\Gamma$ is an open set in
$\Omega$. So $\psi$ is continuous. Since $\Omega, \partial\Gamma$ are compact Hausdorff
spaces, $\psi$ is a homeomorphism. Again, the two actions on $\Omega$ and $\partial\Gamma$ are
defined by left multiplication and hence are conjugate. \hfill 
$\Box$

\medskip

{\bf Remark \quad} Since the action of $\Gamma$ on $\partial\Gamma$ depends only on the
group structure of $\Gamma$ in ~\cite{gph}, the above proposition shows that
$\mathcal{O}_{\Gamma}$ is, up to isomorophism, independent of the choice of
generators of $\Gamma$.

\section{Nuclearity, simplicity and pure infiniteness of $\mathcal{O}_{\Gamma}$}
We first begin by reviewing the crossed product $B\rtimes\mathbb{N}$ of a $C^\ast$-algebra 
$B$ by a $\ast$-endomorphism; this construction was first introduced by Cuntz \cite{cun} to describe the Cuntz algebra
$\mathcal{O}_n$ as the crossed product of UHF algebras by
$\ast$-endomorphisms. See Stacey's paper \cite{sta} for a more detailed
discussion. Suppose that $\rho$ is an injective $\ast$-endomorphism on a unital 
$C^\ast$-algebra $B$. Let $\overline{B}$ be the inductive limit $\varinjlim(B\stackrel{\rho}
{\longrightarrow}B)$ with the corresponding injective homomorphisms $\sigma_n : 
B\to\overline{B}$ $(n\in\mathbb{N})$. Let $p$ be the projection $\sigma_0(1)$. There exists an 
automorphism $\bar{\rho}$ given by $\bar\rho\circ\sigma_n=\sigma_n\circ\rho$ with inverse 
$\sigma_n(b)\mapsto\sigma_{n+1}(b)$. Then the crossed product $B\rtimes_{\rho}\mathbb{N}$ is 
defined to be the hereditary $C^\ast$-algebra 
$p(\overline{B}\rtimes_{\bar\rho}\mathbb{Z})p$. The map $\sigma_0$ induces an embedding of 
$B$ into $\overline{B}$. Therefore the canonical embedding of $\overline{B}$ into 
$\overline{B}\rtimes_{\bar\rho}\mathbb{Z}$ gives an embedding $\pi : B\to 
B\rtimes_{\rho}\mathbb{N}$. Moreover the compression by $p$ of the implementing unitary is 
an isometry $V$ belonging to $B\rtimes_{\rho}\mathbb{N}$ satisfying 
$$V\pi(b)V^\ast=\pi(\rho(b)).$$
In fact, $B\rtimes_{\rho}\mathbb{N}$ is also the universal $C^\ast$-algebra generated by a 
copy $\pi(B)$ of $B$ and an isometry $V$ satisfying the above relation. If $B$ is 
nuclear, then so is $B\rtimes_{\rho}\mathbb{N}$.

\begin{proposition}
$$\mathcal{O}_{\Gamma}\simeq\mathcal{O}_{\Gamma}^{\mathbb{T}}\rtimes_{\rho}\mathbb{N}$$

In particular, $\mathcal{O}_{\Gamma}$ is nuclear.
\end{proposition}
{\it Proof.}
We fix $g_i\in G_i\setminus H$ for all $i\in I$. We
can choose projections $e_i$ which are sums of projections $P_g$ such that $e_i\leq
Q_{g_i}$ and $\sum_{i\in I}e_i=1$. Then $V=\sum_{i\in 
I}S_{g_i}e_i$ is an isometry in
$\mathcal{O}_{\Gamma}$. 

We claim that
$V\mathcal{O}_{\Gamma}^{\mathbb{T}}V^{\ast}\subseteq\mathcal{O}_{\Gamma}^{\mathbb{T}}$
and
$\mathcal{O}_{\Gamma}=C^\ast\left(\mathcal{O}_{\Gamma}^{\mathbb{T}},
  V\right)$. Let $a\in\mathcal{O}_{\Gamma}^{\mathbb{T}}$. It is
obvious that $VaV^\ast\in\mathcal{O}_{\Gamma}^{\mathbb{T}}$ and $C^\ast\left(\mathcal{O}_{\Gamma}^{\mathbb{T}},
  V\right)\subseteq\mathcal{O}_{\Gamma}$. To show the second claim, it suffices to
check that $S_{\mu}P_iS_{\nu}^\ast\in\mathcal{O}_{\Gamma}$ for all $\mu, \nu$
and $i$. If $|\mu|=|\nu|$, we have $S_{\mu}P_iS_{\nu}^\ast\in\mathcal{O}_{\Gamma}^{\mathbb{T}}$.
If $|\mu|\ne|\nu|$, then we may
assume $|\mu|<|\nu|$. Let $|\nu|-|\mu|=k$. Thus $S_{\mu}P_iS_{\nu}^\ast=(V^\ast)^kV^kS_{\mu}P_iS_{\nu}^\ast$ 
and
$V^kS_{\mu}P_iS_{\nu}^\ast\in\mathcal{O}_{\Gamma}^{\mathbb{T}}$.
This proves our claim.

We define a $\ast$-endomorphism $\rho$ of
$\mathcal{O}_{\Gamma}^{\mathbb{T}}$ by $\rho(a)=VaV^\ast$ for
$a\in\mathcal{O}_{\Gamma}^{\mathbb{T}}$. Thanks to the universality of
the crossed product
$\mathcal{O}_{\Gamma}^{\mathbb{T}}\rtimes_{\rho}\mathbb{N}$, we
obtain a canonical surjective $\ast$-homomorphism $\sigma$ of
$\mathcal{O}_{\Gamma}^{\mathbb{T}}\rtimes_{\rho}\mathbb{N}$ onto
$C^\ast(\mathcal{O}_{\Gamma}^{\mathbb{T}}, V)$. Since
$\mathcal{O}_{\Gamma}^{\mathbb{T}}\rtimes_{\rho}\mathbb{N}$ has the
universal property,
there also exists a gauge action $\beta$ on
$\mathcal{O}_{\Gamma}^{\mathbb{T}}\rtimes_{\rho}\mathbb{N}$. Let
$\Psi$ be the corresponding canonical conditional expectation of
$\mathcal{O}_{\Gamma}^{\mathbb{T}}\rtimes_{\rho}\mathbb{N}$ onto
$\mathcal{O}_{\Gamma}^{\mathbb{T}}$. Suppose that $a\in{\rm ker}\sigma$. Then $\sigma(a^\ast 
a)=0$. Since
$\alpha\circ\sigma=\sigma\circ\beta$, we have $\sigma\circ \Psi(a^\ast
a)=0$. The injectivity of $\sigma$ on
$\mathcal{O}_{\Gamma}^{\mathbb{T}}$ implies $\Psi(a^\ast
a)=0$ and hence $a^\ast
a=0$ and $a=0$. It follows that $\mathcal{O}_{\Gamma}\simeq\mathcal{O}_{\Gamma}^{\mathbb{T}}
\rtimes_{\rho}\mathbb{N}$. \hfill $\Box$

\medskip

 In section 2, we reviewed the notion of amenability for discrete group actions. The following is a special case of ~\cite{ada}.

\begin{corollary}
The action of $\Gamma$ on $\partial\Gamma$ is amenable.
\end{corollary}
{\it Proof.}
This follows from Theorem 2.2 and the above proposition. \hfill $\Box$

\medskip
 
We also have a partial result of ~\cite{kir}, ~\cite{dyk}, ~\cite{dyk2} and ~\cite{ds}.

\begin{corollary}
The reduced group $C^\ast$-algebra $C^\ast_r(\Gamma)$ is exact.
\end{corollary}
{\it Proof.}
It is well-known that every $C^\ast$-subalgebra of an exact $C^\ast$-algebra is exact; see 
Wassermann's monograph ~\cite{was}.
Therefore the inclusion
$C_r^\ast(\Gamma)\subset\mathcal{O}_{\Gamma}$ implies exactness. \hfill $\Box$

\medskip

Finally we give a sufficient condition for the 
simplicity and pure infiniteness of $\mathcal{O}_\Gamma$.

\begin{corollary}
Suppose that $\Gamma=\ast_H G_i$ satisfies the following
condition:

There exists at least one element $j\in I$ such that
$$\bigcap_{i\ne j}N_i=\{e\}, $$
where $N_i=\bigcap_{g\in G_i}gHg^{-1}$.

Then $\mathcal{O}_{\Gamma}$ is simple and purely infinite.
\end{corollary}
{\it Proof.}
We first claim that for any $\mu\in\Delta$ and $|g|=1$ with $|\mu g|=|\mu|+1$, 
$$\mu H\mu^{-1}\cap H\supseteq\mu gHg^{-1}\mu^{-1}\cap H.$$
Suppose that $\mu=\mu_1\cdots\mu_n$ such that $\mu_k\in\Omega_{i_k}$
with $\mu_1\ne\dots\ne\mu_n$ and $g\in G_i$ with $i\ne i_n$. We first
assume that $\mu=\mu_1$. If $\mu ghg^{-1}\mu^{-1}\in\mu gHg^{-1}\mu^{-1}\cap H$,
then $ghg^{-1}\in\mu^{-1}H\mu\subseteq G_{i_1}$. Thus $ghg^{-1}\in
G_i\cap G_{i_1}$ implies $ghg^{-1}\in H$. Next we assume that $|\mu|>1$.
If $\mu ghg^{-1}\mu^{-1}\in\mu gHg^{-1}\mu^{-1}\cap H$,
then 
$$\mu_2\cdots\mu_nghg^{-1}\mu_k^{-1}\cdots\mu_2^{-1}\in
\mu_1^{-1}H\mu_1\subseteq G_{i_1}.$$ 
Thus $|\mu_2\cdots\mu_nghg^{-1}\mu_k^{-1}\cdots\mu_2^{-1}|\leq 1$ implies $ghg^{-1}\in H$. This proves the claim.

Let $\{S_g, U_h\}$ be any family satisfying the relations (1), (2),
(3) and (4). By the uniqueness theorem, it is enough to show that 
$C^\ast(P_iU_hP_i \mid h\in H)\simeq C^\ast(H)$ for any $i\in I$. We
next claim that there exists $\nu\in\Gamma$ such that
the initial letter of $\nu$ belongs to $\Omega_i$ and 
$\{U_hS_{\nu}\}_{h\in H}$ have mutually orthogonal ranges. 

Let $g\in\Omega_i$. If $gHg^{-1}\cap H=\{e\}$, then it is enough to 
set $\nu=g$. Now suppose that there exists some $h\in gHg^{-1}\cap H$ with $h\ne
e$. We first assume that $i=j$. By the hypothesis, there exists some
$i_1\in I$ such that $g^{-1}hg\not\in N_{i_1}$ and $i\ne i_1$. Hence
there exists $g_1\in\Omega_{i_1}$ such that $g^{-1}hg\not\in g_1Hg_1^{-1}$ 
and so $h\not\in gg_1Hg_1^{-1}g^{-1}$. If $gg_1Hg_1^{-1}g^{-1}\cap H=\{e\}$, 
then it is enough to put $\nu=gg_1$. If not, we set $\gamma_1=g_1g_1'$ 
for some $g_1'\in\Omega_j$. By the first part of the proof, we
have
$$gHg^{-1}\cap H\supsetneqq\mu\gamma_1H\gamma_1^{-1}\mu^{-1}\cap H.$$
Since $H$ is finite, we can inductively obtain $\gamma_1, \gamma_2,
\dots\gamma_n$ satisfying
$$ gHg^{-1}\cap H \supsetneqq g\gamma_1H\gamma_1^{-1}g^{-1}\cap H
\supsetneqq\cdots\supsetneqq g\gamma_1\cdots\gamma_nH\gamma_n^{-1}\cdots\gamma_1^{-1}g^{-1}\cap
H=\{e\}.$$ 
Then we set $\nu=g\gamma_1\cdots\gamma_n$. If $i\ne j$, we can carry out the same arguments by replacing $g$ by $\gamma=gg_j$ for some $g_j\in
\Omega_j$. Hence from the identification $U_hS_{\nu}\leftrightarrow\delta_h\in
l^2(H)$, it follows that the unitary representation $P_iU_hP_i$ is 
quasi-equivalent to the left regular representation of $H$. Thus $\mathcal{O}_{\Gamma}$ is 
simple.

In Section 5, we have proved that $\mathcal{O}_{\Gamma}\simeq
C(\Omega)\rtimes_r\Gamma$. We show that the action of $\Gamma$ on $\Omega$ is the strong 
boundary
action (see Preliminaries). Let $U, V$ be any non-empty open sets in
$\Omega$. There exists some open set $O=\{(x(n))\in\Omega\mid x(1)=g_1,
\cdots, x(k)=g_k\}$ which is contained in $V$. We may also assume that 
$U^c$ is an open of the form $\{(x(n))\in\Omega\mid x(1)=\gamma_1,
\cdots, x(m)=\gamma_m\}$. Let $\gamma=g_1\cdots
g_k\gamma_m^{-1}\cdots\gamma_1^{-1}$. Then we have $\gamma U^c\subset O\subset
V$. Since $C(\Omega)\rtimes_r\Gamma$ is simple, it follows from
~\cite{as} that the action of $\Gamma$ is topological free. Therefore
it follows from Theorem 2.4 that $C(\Omega)\rtimes_r\Gamma$, namely
$\mathcal{O}_{\Gamma}$, is purely infinite. \hfill $\Box$

\medskip

{\bf Remark \quad} We gave a sufficient condition for
$\mathcal{O}_{\Gamma}$ to be simple. However, we can completely 
determine the ideal structure of $\mathcal{O}_{\Gamma}$ with further effort. Indeed, we will obtain
a matrix $A_{\Gamma}$ to compute K-groups of $\mathcal{O}_{\Gamma}$ in 
the next section. The same argument as in ~\cite{ck2} also works for the ideal 
structure of
$\mathcal{O}_{\Gamma}$. For Cuntz-Krieger algebras, we need to assume
that corresponding matrices have the condition (II) of ~\cite{ck2} to apply the
uniqueness theorem. Since we have another uniqueness
theorem for our algebras, we can always apply the ideal structure theorem. 

Let $\Sigma=I\times \{1, \dots , r\}$ be a finite set, where $r$ is the number of all irreducible unitary representations of $H$. For $x, y\in\Sigma$, we define $x\geq y$ if there exists a sequence $x_1, \dots , x_m$ of elements in $\Sigma$ such that $x_1=x, x_m=y$ and $A_{\Gamma}(x_a, x_{a+1})\ne 0 (a=1, \dots , m-1).$ We call $x$ and $y$ equivalent if $x\geq y\geq x$ and write $\varGamma_{A_\Gamma}$ for the partially ordered set of equivalence classes of elements $x$ in $\Sigma$ for which $x\geq x$. A subset $K$ of $\varGamma_{A_\Gamma}$ is called hereditary if $\gamma_1\geq \gamma_2$ and $\gamma_1\in K$ implies $\gamma_2\in K$. Let 
$$\Sigma(K)=\{ x\in\Sigma \mid x_1\geq x\geq x_2 \quad \mbox{for some} \quad x_1, x_2\in\bigcup_{\gamma\in K}\gamma \}.$$
We denote by $I_K$ the closed ideal of $\mathcal{O}_{\Gamma}$ generated by projections $P(i, k)$, which is defined in the next section, for all $(i, k)\in\Sigma(K)$.

\begin{theorem}[{~\cite[Theorem 2.5.]{ck2}}]
The map $K\mapsto I_K$ is an inclusion preserving bijection of the set of hereditary subsets of $\varGamma_{A_\Gamma}$ onto the set of closed ideals of $\mathcal{O}_\Gamma$.
\end{theorem}

\medskip

\section{$K$-theory for $\mathcal{O}_{\Gamma}$}

In this section we give explicit formulae of the $K$-groups of $\mathcal{O}_{\Gamma}$. We 
have described $\mathcal{O}_{\Gamma}$ as
the crossed product
$\mathcal{O}_{\Gamma}^{\mathbb{T}}\rtimes\mathbb{N}$ in Section 6. So to apply the 
Pimsner-Voiculescu exact sequence ~\cite{pv}, 
we need to compute the $K$-groups of the $AF$-algebra
$\mathcal{O}_{\Gamma}^{\mathbb{T}}$. We assume that each $G_i$ is finite for
simplicity throughout this section. We can also compute the $K$-groups for
general cases by essentially the same arguments. Recall that the fixed-point
algebra is described as  follows:

$$\mathcal{O}_{\Gamma}^{\mathbb{T}}=\overline{\bigcup_{n\geqq
0}\mathcal{F}_n},$$
$$\mathcal{F}_n=\oplus_{i\in I} \mathcal{F}_n^i.$$
For each
$n$, we consider a direct summand of $\mathcal{F}_n$, which is
$$\mathcal{F}_n^i=C^\ast( \,
S_{\mu}P_iU_hP_iS_{\nu}^\ast \mid h\in H, |\mu|=|\nu|=n
\, ),$$
and the embedding $\mathcal{F}_n^i\hookrightarrow\mathcal{F}_{n+1}$ is
given by
\begin{eqnarray*}
   \lefteqn{S_{\mu}P_iU_hP_iS_{\nu}^\ast} \\
      &=& \sum_{g\in\Omega_i\setminus\{e\}}S_{\mu}U_h(S_gQ_gS_g^\ast)S_{\nu}^\ast \\
      &=& \sum_g\sum_{i'\ne i}S_{\mu}S_{hg}P_{i'}S_{\nu g}^\ast .
\end{eqnarray*}

Let $\{\chi_1 , \ldots , \chi_r\}$ be the set of characters
corresponding with all irreducible unitary representations of the finite
group $H$ with degrees $n_1 ,\ldots , n_r$. Then we have the identification $C^\ast(H)\simeq
M_{n_1}(\mathbb{C})\oplus\cdots\oplus M_{n_r}(\mathbb{C})$. We can
write a unit $p_k$ of the $k$-th component $M_{n_k}(\mathbb{C})$ of
$C^\ast(H)$ as follows:
$$p_k=\frac{n_k}{|H|} \sum_{h\in H} \overline{\chi_k(h)}U_h.$$

Suppose that for $i\ne j$,
$$\mathcal{F}_n^i\simeq M_{N(n, i)}(\mathbb{C})\otimes C^\ast(H),$$
$$\mathcal{F}_{n+1}^j\simeq M_{N(n+1, j)}(\mathbb{C})\otimes C^\ast(H).$$
Now we compute each embedding of 
$\mathcal{F}_n^i\hookrightarrow\mathcal{F}_{n+1}^j$,
$$M_{N(n, i)}(\mathbb{C})\otimes M_{n_i}(\mathbb{C})\hookrightarrow M_{N(n+1, j)}(\mathbb{C})\otimes 
M_{n_j}(\mathbb{C})$$ 
at the $K$-theory level. $P(i, k)$ denotes $P_ip_kP_i$. Let $P$ be the projection $e\otimes 1$ in $M_{N(n, i)}(\mathbb{C})\otimes 
M_{n_k}(\mathbb{C})$ given by
$$P=S_{\mu}P(i, k)S_{\mu}^\ast \qquad \mbox{for some}\ 
\mu\in\Delta_n,$$
where $e$ is a minimal projection in the matrix algebras, and $Q$ be the unit of $M_{N(n+1, j)}(\mathbb{C})\otimes M_{n_l}(\mathbb{C})$
given by
$$Q=\sum_{\nu\in\Delta_{n+1}}S_{\nu}P(j, l)S_{\nu}^\ast.$$
At the $K$-theory level, we have $[P]=n_k[e]$. Hence it suffices to compute
${\rm tr}(PQ)/n_k$, where $\rm tr$ is
the canonical trace in the matrix algebras.

\begin{eqnarray*}
\frac{{\rm tr}(PQ)}{n_k} &=& {\rm tr} \left( 
\frac{1}{n_k}(S_{\mu}P(i,
k)S_{\mu}^\ast)(\sum_{\nu\in\Delta_{n+1}}S_{\nu}P(j, l)S_{\nu}^\ast) \right) \\
 &=& {\rm tr} \left( \frac{1}{|H|}(\sum_{h\in H}\overline{\chi_k(h)}
 S_{\mu}U_hP_iS_{\mu}^\ast)(\sum_{\nu\in\Delta_{n+1}}S_{\nu}P(j, l)S_{\nu}^{\ast} ) \right) \\
 &=& \frac{1}{|H|} {\rm tr} \left( \sum_{h\in H}\overline{\chi_k(h)}(
 \sum_{g\in\Omega_i\setminus\{e\}}\sum_{i'\ne i}S_{\mu}S_{hg}P_{i'}S_{\mu
 g}^\ast)(\sum_{\nu\in\Delta_{n+1}}S_{\nu}P(j, l)S_{\nu}^{\ast} ) \right) \\
 &=& \frac{1}{|H|} {\rm tr} \left( \sum_{h\in H}\overline{\chi_k(h)}(
 \sum_{g\in\Omega_i\setminus\{e\}}S_{\mu}S_{hg}P(j, l)S_{\mu g}^\ast )\right) 
\\
 &=& \frac{1}{|H|} \sum_{g\in\Omega_i\setminus\{e\}}\sum_{h\in
 H(g)}\overline{\chi_k(h)} {\rm tr} \left(S_{\mu g}U_{g^{-1}hg}P(j, l)S_{\mu
 g}^\ast\right) \\
  &=& \frac{1}{|H|}\sum_{g\in\Omega_i\setminus\{e\}}\sum_{h\in
 H(g)}\overline{\chi_k(h)}\chi_l({g}^{-1}hg),
\end{eqnarray*}
where $H(g)$ is the stabilizer of $gH$ by the left multiplication of
$H$. 

Now fix $x\in X_i\setminus\{e\}$. Let $\{g\in\Omega_i \mid HgH=HxH\}=\{g_0=x, 
g_1,
\ldots , g_{m-1}\}$. Then there exists $h_1, h'_1, \ldots , h_{m-1},
h'_{m-1} \in H$ such that $h_1x=g_1h'_1, \ldots$,
$h_{m-1}x=g_{m-1}h'_{m-1}$. Note that $h_s H(x) h_s^{-1}=H(g_s)$ for
$s=1, \ldots , m-1$. Since $\chi_k, \chi_l$ are class
functions, we have
\begin{eqnarray*}
   \frac{{\rm tr}(PQ)}{n_k} &=& \frac{1}{|H|}\sum_{x\in 
X_i}\left(\sum_{s=1}^{m-1}\sum_{h\in 
H(x)}\overline{\chi_k(h_shh_s^{-1})}\chi_l(h'_sx^{-1}h_s^{-1}\cdot
h_shh_s^{-1}\cdot h_sx{h'}_s^{-1})\right) \\
      &=& \frac{1}{|H|}\sum_{x\in X_i}\left(\sum_{s=1}^{m-1}\sum_{h\in 
H(x)}\overline{\chi_k(h_shh_s^{-1})}\chi_l(h'_sx^{-1}hx{h'}_s^{-1})\right) \\
      &=& \frac{1}{|H|}\sum_{x\in X_i}\left(\sum_{s=1}^{m-1}\sum_{h\in 
H(x)}\overline{\chi_k(h)}\chi_l(x^{-1}hx)\right) \\
      &=& \frac{1}{|H|}\sum_{x\in
      X_i}\left(\sum_{s=1}^{m-1}\sum_{h\in
      H(x)}\overline{\chi_k(h)}\chi_l^x(h)\right) \\
      &=& \sum_{x\in
      X_i}\left(\frac{|H(x)|}{|H|}\sum_{s=1}^{m-1}\langle\chi_k, 
\chi_l^x\rangle_{H(x)}\right) \\
      &=& \sum_{x\in X_i}\langle\chi_k, \chi_l^x\rangle_{H(x)},
\end{eqnarray*}
where 
\begin{eqnarray*}
\chi_l^x(h) &=& \chi_l\left(x^{-1}hx\right) \\
\langle\chi_k, \chi_l^x\rangle_{H(x)}&=& \frac{1}{|H(x)|}\sum_{h\in H(x)}
\overline{\chi_k(h)}\chi_l^x(h).
\end{eqnarray*}

Let $A_{\Gamma}((j, l), (i, k))=\sum_{x\in
  X_i\setminus\{e\}}\langle\chi_k, \chi_l^x\rangle_{H(x)}$ for
  $i\ne j$ and $A_{\Gamma}((i, k), (i, l))=0$ for $1\leq k, l \leq r$. Then we describe the embedding
$\mathcal{F}_n^{i}\hookrightarrow\mathcal{F}_{n+1}^j$
 at the $K$-theory level by the matrix $[A_{\Gamma}((i, k), (j, l))]_{1\leq k,
  l\leq r}$. Let
 $A_{\Gamma}=[A_{\Gamma}((i, k), (j, l))]$. We have the following lemma.

\begin{lemma}
$$K_0\left(\mathcal{O}_\Gamma^{\mathbb{T}}\right)=\varinjlim\left(
 \, \mathbb{Z}^N\stackrel{A_{\Gamma}}{\longrightarrow}\mathbb{Z}^N \, \right)$$
$$K_1\left( \, \mathcal{O}_\Gamma^{\mathbb{T}} \, \right)=0$$
where $N=|I|r$.
\end{lemma}

We can compute the $K$-groups of $\mathcal{O}_\Gamma$ by using the
Pimsner-Voiculescu sequence with essentially the same argument as in the
Cuntz-Krieger algebra case (see ~\cite{ck2}).

\begin{theorem}
\begin{eqnarray*}
K_0(\mathcal{O}_\Gamma) &=& \mathbb{Z}^N/(1-A_{\Gamma})\mathbb{Z}^N. \\
K_1(\mathcal{O}_\Gamma) &=& {\rm Ker}\{ 1-A_{\Gamma}:\mathbb{Z}^N\to\mathbb{Z}^N\} \quad \mbox{on}\ \mathbb{Z}^N.
\end{eqnarray*}
\end{theorem}
{\it Proof.} 
It suffices to 
 compute the $K$-groups of 
$\overline{\mathcal{O}}_{\gamma}=\overline{\mathcal{O}}_{\Gamma}^{\mathbb{T}}\rtimes_{\bar{\rho}}\mathbb{Z}$. We represent 
 the inductive limit 
$$\varinjlim\left(
 \, \mathbb{Z}^N\stackrel{A_{\Gamma}}{\longrightarrow}\mathbb{Z}^N \, \right)$$
as the set of equivalence classes of $x=(x_1, x_2, \cdots)$ such that
$x_k\in\mathbb{Z}^{N}$ with $x_{k+1}=A(x_k)$. If $S$ is a partial
 isometry in $\mathcal{O}_{\Gamma}$ such that
 $\alpha_z(S)=zS$ and $P$ is a projection in
 $\mathcal{O}_{\Gamma}^\mathbb{T}$ with $P\leq S^\ast S$,
 then $[\rho(P)]=[VPV^\ast]=[(VS^\ast S)P(VS^\ast S)^\ast]=[SPS^\ast]$ in
 $K_0(\mathcal{O}_{\Gamma}^\mathbb{T})$. Recall that 
$$p_k=\frac{n_k}{|H|} \sum_{h\in H} \overline{\chi_k(h)}U_h.$$
Let $P=S_{\mu}P(i, k)S_{\mu}^\ast$  for some $\mu\in\Delta_n$.
If $\mu=\mu_1\cdots\mu_n$, then
\begin{eqnarray*}
   \lefteqn{[\bar{\rho}^{-1}(P)]} \\
   &=&[S_{\mu_1}^\ast PS_{\mu_1}] \\
   &=&[\frac{n_k}{|H|} \sum_{h\in H} \overline{\chi_k(h)}
   \left(S_{\mu_2}\cdots
   S_{\mu_n}P_iU_hP_iS_{\mu_n}\cdots
   S_{\mu_2}^\ast\right)] \\
   &=&\cdots \\
   &=&\sum_{j\ne i}\sum_{l=1}^{r}n_i\left(\sum_{x\in X_i\setminus\{e\}}\langle\chi_k, \chi_l^x\rangle[e_l]\right), 
\end{eqnarray*} 
where the $e_l$ are non-zero minimal projections for $1\leq l \leq
r$. Thus it follows that $\bar{\rho}_\ast^{-1}$ is the shift on $K_0(\overline{\mathcal{O}}_{\Gamma}^{\mathbb{T}})$. We denote the shift
by $\sigma$. If $x=(x_1, x_2, x_3, \cdots )\in
K_0(\overline{\mathcal{O}}_{\Gamma}^{\mathbb{T}})$, then
$\sigma(x)=(x_2, x_3, \cdots )$. By the
Pimsner-Voiculescu exact sequence, there exists an exact sequence
$$0\to K_1(\overline{\mathcal{O}}_{\Gamma})\to
K_0(\overline{\mathcal{O}}_{\Gamma}^{\mathbb{T}})\to
K_0(\overline{\mathcal{O}}_{\Gamma}^{\mathbb{T}})\to
K_0(\overline{\mathcal{O}}_{\Gamma})\to 0. $$
It therefore follows that
$K_0(\overline{\mathcal{O}}_{\Gamma})=K_0(\overline{\mathcal{O}}_{\Gamma}^{\mathbb{T}})/(1-\sigma)K_0(\overline{\mathcal{O}}_{\Gamma}^{\mathbb{T}})$ and
$K_1(\overline{\mathcal{O}}_{\Gamma})={\rm ker}(1-\sigma)$ on
$K_0(\overline{\mathcal{O}}_{\Gamma}^{\mathbb{T}})$. \hfill $\Box$

\medskip

Finally we consider some simple examples.
First let $\Gamma=SL(2,
\mathbb{Z})=\mathbb{Z}_4\ast_{\mathbb{Z}_2}\mathbb{Z}_6$. Let $\chi_1$ 
be the unit character of $\mathbb{Z}_2$ and let $\chi_2$ be the
character such that $\chi_2(a)=-1$ where $a$ is a generator of
$\mathbb{Z}_2$. These are one-dimensional and exhaust all the irreducible characters. 
Then we have the corresponding matrix
$$A_{\Gamma}=\left(\begin{array}{cccc}
0 & 0 & 1 & 0 \\
0 & 0 & 0 & 1 \\
2 & 0 & 0 & 0 \\
0 & 2 & 0 & 0
\end{array}\right) .$$
Hence the corresponding $K$-groups are $K_0(\mathcal{O}_{\Gamma})=0$ and 
$K_1(\mathcal{O}_{\Gamma})=0$. In fact, 
$\mathcal{O}_{\mathbb{Z}_4\ast_{\mathbb{Z}_2}\mathbb{Z}_6}\simeq\mathcal{O}_{\mathbb{Z}_2\ast\mathbb{Z}_3}\oplus\mathcal{O}_{\mathbb{Z}_2\ast\mathbb{Z}_3}\simeq\mathcal{O}_2\oplus\mathcal{O}_2$.

Next let
$\Gamma=\mathfrak{S}_4\ast_{{\mathfrak{S}}_3}\mathfrak{S}_4$, $\tau=(1\, 2)$ 
and $\sigma=(1 \, 2 \, 3)$. Note that $\mathfrak{S}_3=\langle 1, \tau,
\sigma\rangle$. ${\mathfrak{S}}_3$ has three irreducible characters:

$$\tabcolsep=3mm
\begin{tabular}{|c|c|c|c|} \hline
       & 1 & $\tau$ & $\sigma$ \\ \hline
$\chi_1$ & 1 & 1    & 1      \\ \hline
$\chi_2$ & 1 & $-1$    & 1      \\ \hline
$\chi_3$ & 2 & 0     & $-1$ \\ \hline
\end{tabular}$$

Moreover, $\mathfrak{S}_3\backslash\mathfrak{S}_4/\mathfrak{S}_3$ has
only two points; say $\mathfrak{S}_3$ and $\mathfrak{S}_3 x\mathfrak{S}_3$ with $x=(1\, 
2)(3\, 4)$. Then we obtain the corresponding matrix
$$A_{\Gamma}=\left(\begin{array}{cccccc}
0 & 0 & 0 & 1 & 0 & 1 \\
0 & 0 & 0 & 0 & 1 & 1 \\
0 & 0 & 0 & 1 & 1 & 2 \\
1 & 0 & 1 & 0 & 0 & 0 \\
0 & 1 & 1 & 0 & 0 & 0 \\
1 & 1 & 2 & 0 & 0 & 0
\end{array}\right) .$$
Hence this gives
$K_0(\mathcal{O}_{\Gamma})=\mathbb{Z}\oplus\mathbb{Z}_4$ and
$K_1(\mathcal{O}_{\Gamma})=\mathbb{Z}$. In this case, $\Gamma$
satisfies the condition of Theorem 6.3. So $\mathcal{O}_{\Gamma}$ is a
simple, nuclear, purely infinite $C^\ast$-algebra.

\section{KMS states on $\mathcal{O}_{\Gamma}$}

In this section, we investigate the relationship between KMS states on
$\mathcal{O}_{\Gamma}$ for generalized gauge actions and random walks on $\Gamma$. Throughout this section, we assume that all groups $G_i$
are finite though we can carry out the same arguments if $G_i=\mathbb{Z}\times H$ for some $i\in I$. Let $\omega=(\omega_i)_{i\in
  I}\in\mathbb{R}_{+}^{|I|}$. By the universality of $\mathcal{O}_{\Gamma}$, we can define an automorphism $\alpha_t^\omega$ 
for any
$t\in\mathbb{R}$ on $\mathcal{O}_{\Gamma}$ by
$\alpha_t^{\omega}(S_g)=e^{\sqrt{-1} \, \omega_i t}S_g$ for $g\in
G_i\setminus H$ and $\alpha_t^{\omega}(U_h)=U_h$ for $h\in H$. Hence we
obtain the $\mathbb{R}$-action $\alpha^\omega$ on
$\mathcal{O}_{\Gamma}$. We call it {\it the generalized gauge action}
with respect to $\omega$. We will only consider actions of these
types and determine KMS states on
$\mathcal{O}_{\Gamma}$ for these actions.

In ~\cite{woe}, Woess showed that our boundary $\Omega$ can be identified with
the Poisson boundary of random walks satisfying certain conditions. The reader is referred to ~\cite{woe2} for a good survey of random walks.

Let $\mu$ be a probability measure on
$\Gamma$ and consider a random walk governed by $\mu$, i.e. the transition
probability from $x$ to $y$ given by 
$$p(x, y)=\mu(x^{-1}y).$$ 
A random walk is said to be {\it irreducible} if for any $x, y\in\Gamma$,
$p^{(n)}(x, y)\ne 0$ for some integer $n$, where 
$$p^{(n)}(x, y)=\sum_{x_1, x_2, \dots , x_{n-1}\in\Gamma}p(x, x_1)p(x_1,
x_2)\cdots p(x_{n-1}, y).$$
A probability measure $\nu$ on $\Omega$ is said to be {\it stationary}
with respect to $\mu$ if $\nu=\mu\ast\nu$, where $\mu\ast\nu$ is
defined by
$$\int_{\Omega}f(\omega)d\mu\ast\nu(\omega)=\int_{\Omega}\int_{{\rm
 supp}\mu}f(g\omega)d\mu(g)d\nu(\omega),\quad\mbox{for}\quad f\in C(\Omega, \nu).$$
 By {~\cite[Theorem 9.1]{woe}}, if a random walk governed by a probability measure $\mu$ on $\Gamma$ is irreducible, then there exists a unique stationary probability measure $\nu$ on $\Omega$ with
respect to $\mu$. Moreover if $\mu$ has finite support, then the
Poisson boundary coincides with $(\Omega, \nu)$.  

If $\nu$ is a probability measure on the compact
space $\Omega$, then we can define a state $\phi_{\nu}$ by 
$$\begin{array}{ccc}
\phi_{\nu}(X)=\int_{\Omega}E(X) d\nu & \mbox{for} &
 X\in\mathcal{O}_{\Gamma},
\end{array} $$
where $E$ is the canonical conditional expectation of
$C(\Omega)\rtimes_r\Gamma$ onto $C(\Omega)$. 

One of our purposes in this section is to prove that there exists a random walk
governed by a probability measure $\mu$ that induces the stationary measure 
$\nu$ on $\Omega$ such that the corresponding state $\phi_{\nu}$ is
the unique KMS state for $\alpha^{\omega}$. Namely,

\begin{theorem}
Assume that the matrix $A_{\Gamma}$ obtained in the preceding section is irreducible.
For any $\omega=(\omega_i)_{i\in I}\in\mathbb{R}_+^{|I|}$, there
exists a unique probability measure $\mu$ with the following properties:

{\rm (i)} $supp(\mu)=\bigcup_{i\in I}G_i\setminus H$.

{\rm (ii)} $\mu(gh)=\mu(g)$ for any $g\in\bigcup_{i\in I}G_i\setminus H$ and 
$h\in H$.

{\rm (iii)} The corresponding unique stationary measure $\nu$ on $\Omega$
induces the unique KMS state $\phi_\nu$ for $\alpha^{\omega}$ and the
corresponding inverse temperature $\beta$ is also unique.
\end{theorem}

\medskip

We need the hypothesis of the irreducibility of the matrix $A_\Gamma$ for the uniqueness of the KMS state. Though it is, in general, difficult to check the irreducibility of $A_\Gamma$, by Theorem 6.5, the condition of simplicity of $\mathcal{O}_\Gamma$ in Corollary 6.4 is also a sufficient condition for irreducibility of $A_\Gamma$. To obtain the theorem, we first present two lemmas.

\begin{lemma}
Assume that $\nu$ is a probability measure on $\Omega$. Then the
corresponding state $\phi_{\nu}$ is the KMS state for
$\alpha^{\omega}$ if and only if $\nu$ satisfies the following
conditions:
$$\nu(\Omega(x_1\cdots x_m))=\frac{e^{-\beta\omega_{i_1}}\cdots
  e^{-\beta\omega_{i_{m-1}}}}{[G_{i_m}:H]-1+e^{\beta\omega_{i_m}}},$$
for $x_k\in\Omega_{i_k}$ with $i_1\ne\dots\ne i_m$, where
  $\Omega(x_1\cdots x_m)$ is the cylinder subset of $\Omega$ defined by
$$\Omega(x_1\cdots x_m)=\{(x(n))_{n\geq 1}\in\Omega\mid x(1)=x_1, \dots
,x(m)=x_m \}.$$
\end{lemma}

{\it Proof} $\phi_{\nu}$ is the KMS state for
$\alpha^{\omega}$ if and only if 
$$\phi_\nu(S_\xi P_iU_hS_\eta^\ast\cdot S_\sigma P_jU_kS_\tau^\ast)=\phi(S_\sigma P_jU_kS_\tau^\ast\cdot\alpha_{\sqrt{-1}\beta}^\omega(S_\xi P_iU_hS_\eta^\ast)),$$
for any $\xi, \eta, \sigma, \tau\in\Delta, h, k\in H$ and $i, j\in
I$. 

We may assume that $|\xi|+|\sigma|=|\eta|+|\tau|$ and
$|\eta|\geq|\sigma|$. Set $|\xi|=p, |\eta|=q, |\sigma|=s, |\tau|=t$
and let $\xi=\xi_1\cdots\xi_p$, $\eta=\eta_1\cdots\eta_q$ with
$\xi_k\in\Omega_{i_k}\setminus\{e\},
\eta_l\in\Omega_{j_l}\setminus\{e\}$ and $i_1\ne\dots\ne i_p,
j_1\ne\dots\ne j_q$. Then 
\begin{eqnarray*}
\lefteqn{\phi_\nu(S_\xi P_iU_hS_\eta^\ast\cdot S_\sigma
  P_jU_kS_\tau^\ast)=\delta_{\eta_1\cdots\eta_s,
  \sigma}\delta_{\eta_{s+1}, j}\phi_\nu(S_\xi
  P_iU_hS_{\eta_{s+1}\cdots\eta_{q}}^\ast U_kS_\tau^\ast)} \\
&=& \delta_{\eta_1\cdots\eta_s, \sigma}\delta_{\eta_{s+1},
  j}\phi_\nu(S_{\xi h}P_iS_{\tau k^{-1}\eta_{s+1}\cdots\eta_{q}}) \\
&=& \delta_{\eta_1\cdots\eta_s, \sigma}\delta_{\eta_{s+1},
  j}\delta_{\xi h, \tau
  k^{-1}\eta_{s+1}\cdots\eta_{q}}\sum_{x\in\Omega_i\setminus\{e\}}\nu(\Omega(\xi x)),
\end{eqnarray*}
and 
\begin{eqnarray*}
\lefteqn{\phi_\nu(S_\sigma
  P_jU_kS_\tau^\ast\cdot\alpha_{\sqrt{-1}\beta}^\omega(S_\xi
  P_iU_hS_\eta^\ast))}\\
&=& e^{-\beta\omega_{i_1}}\cdots
  e^{-\beta\omega_{i_p}}e^{\beta\omega_{j_1}}\cdots
  e^{\beta\omega_{j_q}}\phi_\nu(S_\sigma P_jU_kS_\tau^\ast\cdot S_\xi
  P_iU_hS_\eta^\ast) \\
&=& e^{-\beta\omega_{i_1}}\cdots
  e^{-\beta\omega_{i_p}}e^{\beta\omega_{j_1}}\cdots
  e^{\beta\omega_{j_q}}\delta_{\tau,
  \xi_1\cdots\xi_t}\delta_{\xi_{t+1}, j}\phi_\nu(S_{\sigma
  k\xi_{t+1}\cdots\xi_{p} h} P_i S_\eta^\ast) \\
&=& e^{-\beta\omega_{i_1}}\cdots
  e^{-\beta\omega_{i_p}}e^{\beta\omega_{j_1}}\cdots
  e^{\beta\omega_{j_q}}\delta_{\tau,
  \xi_1\cdots\xi_t}\delta_{\xi_{t+1}, j}\delta_{\sigma
  k\xi_{t+1}\cdots\xi_{p}h, \eta}\sum_{x\in\Omega_i\setminus\{e\}}\nu(\Omega(\eta x)),
\end{eqnarray*}
where $\delta_{g, i}=1$ only if $g\in G_i\setminus H$. Therefore the
corresponding state $\phi_{\nu}$ is the KMS state for
$\alpha^{\omega}$ if and only if $\nu$ satisfies the following
conditions:
$$\nu(\Omega(\xi_1\dots\xi_p x))=e^{-\beta\omega_{i_1}}\cdots
  e^{-\beta\omega_{i_p}}\nu(\Omega(x)),$$
for $x\in\Omega_i\setminus\{e\}$ with $i\ne i_p$.

Now we assume that $\phi_\nu$ is the KMS state for
$\alpha^{\omega}$. Then for $i\in I$,
\begin{eqnarray*}
\lefteqn{\nu(Y_i)=\phi_\nu(P_i)=\sum_{g\in\Omega_i\setminus\{e\}}\phi_\nu(S_gS_g^\ast)}\\
&=&\sum_{g\in\Omega_i\setminus\{e\}}\phi_\nu(S_g^\ast\alpha_{\sqrt{-1}\beta}^\omega(S_g))\\
&=& e^{-\beta\omega_i}\sum_{g\in\Omega_i\setminus\{e\}}\phi_\nu(Q_g)\\
&=& e^{-\beta\omega_i}\sum_{g\in\Omega_i\setminus\{e\}}\phi_\nu(1-P_i)\\
&=& e^{-\beta\omega_i}([G_i:H]-1)(1-\nu(Y_i)).
\end{eqnarray*}
Hence, 
$$\nu(Y_i)=\frac{[G_i:H]-1}{[G_i:H]-1+e^{\beta\omega_i}}.$$
Moreover, 
\begin{eqnarray*}
\lefteqn{\nu(\Omega(x_1\dots x_m))=\phi_\nu(S_{x_1}\cdots
  S_{x_m}S_{x_m}^\ast\cdots S_{x_1}^\ast)}\\
&=&\phi_\nu(S_{x_m}^\ast\cdots
  S_{x_1}^\ast\alpha_{\sqrt{-1}\beta}^\omega(S_{x_1}\cdots S_{x_m}))\\
&=&e^{-\beta\omega_{i_1}}\cdots
  e^{-\beta\omega_{i_m}}\phi_\nu(Q_{x_m})\\
&=&e^{-\beta\omega_{i_1}}\cdots
  e^{-\beta\omega_{i_m}}(1-\nu(\Omega(Y_{i_m})))\\
&=&\frac{e^{-\beta\omega_{i_1}}\cdots
  e^{-\beta\omega_{i_{m-1}}}}{[G_{i_m}:H]-1+e^{\beta\omega_{i_m}}}.
\end{eqnarray*}
Conversely, suppose that a probability measure $\nu$ satisfies the
condition of this lemma. By the first part of this proof, $\phi_\nu$
is the KMS state for $\alpha^\omega$. \hfill $\Box$

\medskip

\begin{lemma}
Assume that $\nu$ is the unique stationary measure on $\Omega$ with respect to
a random walk on $\Gamma$, governed by a probability measure $\mu$ 
with the conditions {\rm (i), (ii)} in Theorem 8.1. Then $\phi_{\nu}$
is a $\beta$-KMS state for $\alpha^\omega$ if and only if $\mu$
satisfies the following conditions:
$$\displaystyle{\mu(g)=\frac{\prod_{j\ne i}C_j}{\sum_{k\in
      I}(g_k\prod_{l\ne k}C_l)}}\qquad \mbox{for}\quad g\in
G_i\setminus H\quad\mbox{and}\quad i\in I,$$ 
where $g_i=|G_i\setminus H|$ and
      $C_i=(1-e^{-\beta\omega_i})g_i-(1-e^{\beta\omega_i})|H|$ for
      $i\in I$.
\end{lemma} 
{\it Proof}
Assume that $\phi_{\nu}$ is a $\beta$-KMS state for
$\alpha^{\omega}$. For any $f\in C(\Omega)$,
\begin{eqnarray*}
\lefteqn{\int\!\!\!\!\!\int f(\omega) d\nu(\omega)=\int\!\!\!\!\!\int
  f(\omega) d\mu*\nu(\omega)} \\
&=&\int\!\!\!\!\!\int f(g \omega)
  d\nu(\omega)d\mu(g) \\
&=& \int\!\!\!\!\!\int (\lambda_g^\ast f\lambda_g)(\omega) d\nu(\omega)d\mu(g) \\
&=& \sum_{g\in {\rm supp}(\mu)}\mu(g)\phi_{\nu}(\lambda_g^\ast f\lambda_g) \\
&=& \sum_{g\in {\rm supp}(\mu)}\mu(g)\phi_{\nu}(f\lambda_g \alpha_{\sqrt{-1}\beta}^{\omega}(\lambda_g^\ast)),
\end{eqnarray*} 
where $\mathcal{O}_{\Gamma}\simeq C(\Omega)\rtimes_r\Gamma=C^\ast(f,
      \lambda_{\gamma} \mid f\in C(\Omega), \gamma\in\Gamma).$

Put $f=\chi_{\Omega(x)}=P_x$ for $i\in I$ and
      $x\in\Omega_i\setminus\{e\}$. Since
      $\lambda_g=S_g+\sum_{g'\in\Omega_{i'}\setminus H\cup g^{-1}H}S_{gg'}S_{g'}^\ast+S_{g^{-1}}^\ast$ for $g\in
      G_{i'}\setminus H$ and $i'\in I$, we
      have
$$1=\sum_{gH=xH}\mu(g)e^{\beta\omega_i}+\sum_{g\in G_i\setminus H, gH\ne
  xH}\mu(g)+\sum_{g\in G_j\setminus H, j\ne
  i}\mu(g)e^{-\beta\omega_j}$$
for any $i\in I$ and $x\in\Omega_i\setminus\{e\}$. 
Let $x, y\in\Omega_i\setminus\{e\}$ with $xH\ne yH$. Then
$$1=\sum_{gH=xH}\mu(g)e^{\beta\omega_i}+\sum_{gH\ne
  xH}\mu(g)+\sum_{g\in G_j\setminus H, j\ne
  i}\mu(g)e^{-\beta\omega_j},$$
$$1=\sum_{gH=yH}\mu(g)e^{\beta\omega_i}+\sum_{gH\ne
  yH}\mu(g)+\sum_{g\in G_j\setminus H, j\ne
  i}\mu(g)e^{-\beta\omega_j}.$$
By the above equations, we have $\mu(x)=\mu(y)$, and then it follows from hypothesis ${\rm (ii)}$ in Theorem 8.1 that $\mu(g)=\mu_i$ for any $g\in G_i\setminus H$. Therefore we have
$$1=|H|e^{\beta\omega_i}\mu_i+(g_i-|H|)\mu_i+\sum_{j\ne i}g_je^{-\beta\omega_j}\mu_j,$$
for any $i\in I$, where $g_i=|G_i\setminus H|$. Thus by considering the above equations for $i$ and $j\in I$,
$$|H|e^{\beta\omega_i}\mu_i-|H|e^{\beta\omega_j}\mu_j+(g_i-|H|)\mu_i-(g_j-|H|)\mu_j+g_je^{-\beta\omega_j}\mu_j-g_ie^{-\beta\omega_i}\mu_i=0.$$
Hence we obtain the equation,
$$(|H|e^{\beta\omega_i}+g_i-|H|-g_ie^{-\beta\omega_i})\mu_i=(|H|e^{\beta\omega_j}+g_j-|H|-g_je^{-\beta\omega_j})\mu_j.$$
Since $\mu(\bigcup_{i\in I}G_i\setminus
  H)=1$, we have 
$$g_i\mu_i+\sum_{j\ne i}g_j\frac{(1-e^{-\beta\omega_i})g_i-(1-e^{-\beta\omega_i})|H|}{(1-e^{-\beta\omega_j})g_j-(1-e^{-\beta\omega_j})|H|}\mu_i=1.$$
We put $C_i=(1-e^{-\beta\omega_i})g_i-(1-e^{-\beta\omega_i})|H|$ and then
$$(g_i+C_i\sum_{j\ne i}\frac{g_j}{C_j})\mu_i=1.$$
Therefore 
\begin{eqnarray*}
\lefteqn{\mu_i=\frac{1}{g_i+C_i\sum_{j\ne i}g_j/C_j}} \\
&=& \frac{\prod_{j\ne i}C_j}{g_i\prod_{j\ne i}C_j+\sum_{j\ne i}(g_jC_i\prod_{k\ne i, j}C_k)} \\
&=& \frac{\prod_{j\ne i}C_j}{\sum_{k\in I}g_k\prod_{l\ne k}C_l}.
\end{eqnarray*}

On the other hand, let $\nu$ be the probability measure on $\Omega$ satisfying the condition in Lemma 8.2. Then the corresponding state $\phi_\nu$ is the KMS state. It is enough to check that $\mu*\nu=\nu$ by ~\cite{woe}. Since 
$$\nu(\Omega(x_1\cdots x_n))=e^{-\beta\omega_{i_1}}\cdots e^{-\beta\omega_{i_{n-1}}}\nu(\Omega(x_n)),$$
for $x_k\in\Omega_{i_k}\setminus\{e\}$ with $i_1\ne\dots\ne i_n$, we have
\begin{eqnarray*}
\lefteqn{\mu*\nu(\Omega(x_1\cdots x_n))=\int\!\!\!\!\!\int
  \chi_{\Omega(x_1\cdots x_n)}(\omega) d\mu*\nu(\omega)} \\
&=& \sum_{g\in {\rm supp}\mu}\mu(g)\int(\lambda_g^\ast\chi_{\Omega(x_1\cdots x_n)}\lambda_g)(\omega)d\nu(\omega) \\
&=& \sum_{g\in G_{i_1}\setminus H, x_1H=gH}\mu_{i_1}\phi_\nu(S_{x_2}\cdots S_{x_n}S_{x_n}^\ast\cdots S_{x_2}^\ast) \\
&& +\sum_{g\in G_{i_1}\setminus H, x_1H\ne gH}\mu_{i_1}\phi_\nu(S_{g^{-1}x_1}S_{x_2}\cdots S_{x_n}S_{x_n}^\ast\cdots S_{x_2}^\ast S_{g^{-1}x_1}^\ast) \\
&& +\sum_{g\in G_i\setminus H, i\ne i_1}\mu_i\phi_\nu(S_{g^{-1}}S_{x_1}S_{x_2}\cdots S_{x_n}S_{x_n}^\ast\cdots S_{x_2}^\ast S_{x_1}^\ast S_{g^{^1}}^\ast) \\
&=& \left(|H|e^{\beta\omega_{i_1}}\mu_{i_1}+(g_{i_1}-|H|)\mu_{i_1}+\sum_{i\ne i_1}g_ie^{-\beta\omega_i}\mu_i\right)\nu(\Omega(x_1\cdots x_n)) \\
&=& \nu(\Omega(x_1\dots x_n)).
\end{eqnarray*}
\hfill $\Box$

\medskip

To prove the uniqueness of KMS states of $\mathcal{O}_{\Gamma}$, we
need the irreducibility of the matrix $A_{\Gamma}$ (See ~\cite{efw2}
for KMS states on Cuntz-Krieger algebras). Set an irreducible matrix
$B=[B((i,k),(j,l))]=[e^{-\beta\omega_i}A_{\Gamma}^t((i,k),(j,l))]$. Let
$K_{\beta}$ be the set of all $\beta$-KMS states for the action
$\alpha^\omega$. We put
$$L_{\beta}=\{y=[y(i,k)]\in \mathbb{R}^N \mid By=y,\quad y(i,k)\geq 0, 
\quad
\sum_{i\in I}\sum_{k=1}^r n_k y(i,k)=1 \}.$$
We now have the necessary ingredients for the proof of Theorem 8.1. 

\medskip

{\it Proof of Theorem 8.1}
We first prove the uniqueness of the corresponding inverse
temperature. Let $\phi$ be a $\beta$-KMS state for $\alpha^{\omega}$. For $i\in I$,
\begin{eqnarray*}
\lefteqn{\phi(P_i)=\sum_{g\in\Omega_i\setminus\{e\}}\phi(S_gS_g^\ast)
=\sum_{g\in \Omega_i\setminus\{e\}}\phi(S_g^\ast
\alpha_{\sqrt{-1}\beta}^{\omega}(S_g))} \\
&=& e^{-\beta\omega_i}\sum_{g\in \Omega_i\setminus\{e\}}\phi(Q_g) \\
&=& e^{-\beta\omega_i}([G_i:H]-1)(1-\phi(P_i)).
\end{eqnarray*}
Thus $\phi(P_i)=\lambda_i(\beta)/(1+\lambda_i(\beta))$, where
$\lambda_i(\beta)=e^{-\beta\omega_i}([G_i:H]-1)$. Since $\sum_{i\in I}P_i=1$, 
$$|I|-1=\sum_{i\in I}\frac{1}{1+\lambda_i(\beta)}.$$
The function $\sum_{i\in I}1/(1+\lambda_i(\beta))$ is a monotone 
increasing continuous function such that
$$\sum_{i\in I}\frac{1}{1+\lambda_i(\beta)}=\left\{\begin{array}{cc}
\sum_{i\in I}1/[G_i:H] & \mbox{if} \,\,\, \beta=0, \\
|I| & \mbox{if} \,\,\, \beta\to\infty.
\end{array}\right. $$
Since $\sum_{i\in I}1/[G_i:H]\leq |I|/2\leq |I|-1$,
there exists a unique $\beta$ satisfying
$$|I|-1=\sum_{i\in I}\frac{1}{([G_i:H]-1)e^{-\beta\omega_i}+1}.$$
Therefore we obtain the uniqueness of the inverse temperature $\beta$. 

We will next show the uniqueness of the KMS state $\phi_{\nu}$. We claim that $K_{\beta}$ is in one-to-one correspondence with
$\L_{\beta}$. In fact, we define a map $f$ from $K_{\beta}$ to
$\L_{\beta}$ by
$$f(\phi)=[\phi(P(i,k))/n_k].$$
Indeed,
\begin{eqnarray*}
\lefteqn{e^{\beta\omega_i}\phi(P(i,k))=\sum_{g\in\Omega_i\setminus\{e\}}\phi(p_kS_g\alpha_{\sqrt{-1}\beta}^{\omega}(S_g^\ast))} \\
&=& \sum_{g\in\Omega_i\setminus\{e\}}\phi(S_g^\ast p_k S_g) \\
&=& \frac{n_k}{|H|}\sum_{g\in\Omega_i\setminus\{e\}}\sum_{h\in
  H}\overline{\chi_k(h)}\phi(S_g^\ast U_h S_g) \\
&=& \frac{n_k}{|H|}\sum_{g\in\Omega_i\setminus\{e\}}\sum_{h\in
  H(g)}\overline{\chi_k(h)}\phi(Q_gU_{g^{-1}hg}) \\
&=& \frac{n_k}{|H|}\sum_{g\in\Omega_i\setminus\{e\}}\sum_{h\in
  H(g)}\overline{\chi_k(h)}\sum_{j\ne i}\phi(P_jU_{g^{-1}hg}P_j) \\
&=& \frac{n_k}{|H|}\sum_{g\in\Omega_i\setminus\{e\}}\sum_{h\in
  H(g)}\overline{\chi_k(h)}\sum_{j\ne i}\sum_{l=1}^{r}\phi(P(j,
l)U_{g^{-1}hg}P(j,l)).
\end{eqnarray*}
Since $\phi$ is a trace on $C^\ast(P(j, l)U_hP(j, l) \mid h\in H)\simeq
M_{n_l}(\mathbb{C})$ and $M_{n_l}(\mathbb{C})$ has a unique tracial state, we have
$$\phi(P(j, l)U_{g^{-1}hg}P(j,l))=\chi_l(g^{-1}hg)\frac{\phi(P(j,
  l))}{n_l}.$$
Therefore, by the same arguments as in the previous section, we obtain
\begin{eqnarray*}
\lefteqn{e^{\beta\omega_i}\phi(P(i,k))} \\
&=& \frac{n_k}{|H|}\sum_{g\in\Omega_i\setminus\{e\}}\sum_{h\in
  H(g)}\overline{\chi_k(h)}\sum_{j\ne i}\sum_{l=1}^{r}\phi(P(j,
l)U_{g^{-1}hg}P(j,l)) \\
&=& n_k\sum_{x\in X_i\setminus\{e\}}\sum_{j\ne i}\sum_{l=1}^{r}\langle\chi_k,
\chi_l^{x}\rangle_{H(x)}\phi(P(j, l))/n_l \\
&=& n_k \sum_{(j, l)}A_{\Gamma}((j, l), (i, k))\phi(P(j, l))/n_l.
\end{eqnarray*}
Hence this is well-defined.

Suppose that $\nu$ is the probability measure in
Lemma 8.2 and $\phi_{\nu}$ is the induced $\beta$-KMS state for
$\alpha^{\omega}$. Set a vector
$y=[y(i,k)=\phi_{\nu}(P(i,k))/n_k]$. Since $y$ is
strictly positive and $B$ is irreducible, $1$ is the eigenvalue which
dominates the absolute value of all eigenvalue of $B$ by the
Perron-Frobenius theorem. It also follows from the Perron-Frobenius
theorem that $L_{\beta}$ has only one element. Hence $f$ is surjective. 

Let $\phi\in K_{\beta}$. For $\xi=\xi_{i_1}\cdots\xi_{i_n}, \eta=\eta_{j_1}\cdots\eta_{j_n}$ with $i_1\ne\dots\ne i_n, j_1\ne\dots\ne j_n$, $h\in H$ and $i\in I$, 
\begin{eqnarray*}
\lefteqn{e^{\beta\omega_{j_1}}\cdots e^{\beta\omega_{j_n}}\phi(S_{\xi}U_hP_iS_{\eta}^\ast)=\phi(S_{\xi}U_hP_i\alpha_{\sqrt{-1}\beta}^{\omega}(S_{\eta}^\ast))} 
\\
&=& \phi(S_{\eta}^\ast S_{\xi}U_hP_i)=\delta_{\xi, \eta}\phi(U_hP_i) \\
&=& \delta_{\xi, \eta}\sum_{k=1}^r\phi(U_hP(i, k))=\delta_{\xi,
  \eta}\sum_{k=1}^r\chi_k(h)\phi(P(i,k))/n_k,
\end{eqnarray*}
because $\phi$ is a trace on $C^\ast( U_hP(i, k) \mid h\in H)\simeq
M_{n_k}(\mathbb{C})$. If
$f(\phi)=f(\psi)$, then the above calculations imply $\phi=\psi$ on
$\mathcal{O}_{\Gamma}^{\mathbb{T}}$. By the KMS condition,
$\phi(b)=0=\psi(b)$ for
$b\notin\mathcal{O}_{\Gamma}^{\mathbb{T}}$. Thus $\phi=\psi$ and $f$
is injective. Therefore $\phi_{\nu}$
is the unique $\beta$-KMS state for $\alpha^{\omega}$.  \hfill $\Box$
  
\medskip

{\bf Remarks and Examples \quad}Let $\nu$ be the corresponding probability measure with
the gauge action $\alpha$. Under the identification $L^{\infty}(\Omega,
\nu)\rtimes_w\Gamma\simeq\pi_{\nu}(\mathcal{O}_{\Gamma})''$,
we can determine the type of the factor by essentially the same
arguments as in ~\cite{efw2}. If $H$ is trivial, then
$\mathcal{O}_{\Gamma}$ is a Cuntz-Krieger algebra for some irreducible 
matrix with $0$-$1$ entries. In this case, we can always apply the
result in ~\cite{efw2}. This fact generalizes
~\cite{rr}. If $H$ is not trivial, then by using the condition of simplicity 
of $\mathcal{O}_\Gamma$ in Corollary 6.4 to check the irreducibility of the
matrix $A_\Gamma$, we can apply Theorem 8.1. In the special case where $G_i=G$
for all $i\in I$, we can easily determine the type of the factor
$\pi_\nu(\mathcal{O}_\Gamma)''$ for the gauge action. The factor
$\pi_\nu(\mathcal{O}_\Gamma)''$ is of type ${\rm III}_\lambda$ where
$\lambda=1/([G:H]-1)^2$ if $|I|=2$ and $\lambda=1/(|I|-1)([G:H]-1)$ if $|I|>2$.
For instance, let
$\Gamma=\mathfrak{S}_4\ast_{{\mathfrak{S}}_3}\mathfrak{S}_4$. We have already
obtained the matrix $A_{\Gamma}$ in section 7, but we can determine that the
factor $L^{\infty}(\Omega, \nu)\rtimes_w\Gamma$ is of type III$_{1/9}$ without
using $A_{\Gamma}$.  

\medskip

We next discuss the converse. Namely any $\mathbb{R}$-actions that have KMS states induced by a probability measure $\mu$ on $\Gamma$
with some conditions is, 
in fact, a generalized gauge action. 

Let $\mu$ be a given probability measure on $\Gamma$ with ${\rm supp}(\mu)=\bigcup_{i\in
  I}G_i\setminus H$. By ~\cite{woe}, there exists an unique probability measure $\nu$ on
$\Omega$ such that $\mu\ast\nu=\nu$. Let $(\pi_\nu, H_{\nu}, x_{\nu})$
be the GNS-representation of $\mathcal{O}_{\Gamma}$ with respect to the state
$\phi_\nu$. We also denote a vector state of $x_{\nu}$ by $\phi_\nu$.
$$\phi_{\nu}(a)=\langle ax_{\nu}, x_\nu \rangle \qquad \mbox{for}
\quad a\in\pi_\nu(\mathcal{O}_{\Gamma})''.$$
Let $\sigma_t^\nu$ be the modular automorphism group of
$\phi_\nu$.

\begin{theorem}
Suppose that $\mu$ is a probability measure on $\Gamma$ such that
${\rm supp}(\mu)=\bigcup_{i\in I}G_i\setminus H$ and $\mu(g)=\mu(hg)$
for any $g\in\bigcup_{i\in I}G_i\setminus H$, $h\in H$. If $\nu$ is 
the corresponding stationary measure with respect to $\mu$, then there exists $\omega_g\in\mathbb{R}_+$ such that 
$$\sigma_t^\nu(\pi_\nu(S_g))=e^{\sqrt{-1}\omega_gt}\pi_\nu(S_g) \quad
\mbox{for}\quad g\in G_i\setminus H, i\in I,$$
and 
$$\sigma_t^\nu(\pi_\nu(U_h))=\pi_\nu(U_h) \quad
\mbox{for}\quad h\in H.$$
\end{theorem}
{\it Proof}
To prove that
$\sigma_t^\nu(\pi_\nu(S_g))=e^{\sqrt{-1}\omega_gt}\pi_\nu(S_g)$, it suffices to show that there exists $\zeta_g\in\mathbb{R}_+$
such that
$$(*)\qquad\phi_\nu(\pi_\nu(S_g)a)=\zeta_g\phi_\nu(a \pi_\nu(S_g))
\quad\mbox{for}\quad g\in G_i\setminus, a\in\pi_\nu(\mathcal{O}_\Gamma)''.$$
In fact, Let $\Delta_\nu$ be the modular operator and $J_\nu$ be the
modular conjugate of $\phi_\nu$.
\begin{eqnarray*}
\lefteqn{\mbox{(left hand side of $(*)$)}=\langle \pi_\nu(S_g)a x_\nu,
  x_\nu\rangle} \\
&=& \langle a x_\nu, \pi_\nu(S_g)^\ast x_\nu\rangle \\
&=& \langle a x_\nu, J_\nu\Delta_\nu^{1/2}\pi_\nu(S_g)x_\nu\rangle \\
&=& \langle \Delta_\nu^{1/2}\pi_\nu(S_g)x_\nu, J_\nu a x_\nu\rangle \\
&=& \langle \Delta_\nu^{1/2}\pi_\nu(S_g)x_\nu, \Delta_\nu^{1/2} a^\ast x_\nu\rangle.
\end{eqnarray*}
and
\begin{eqnarray*}
\lefteqn{\mbox{(right hand side of $(*)$)}=\zeta_g\langle a\pi_\nu(S_g)x_\nu, x_\nu\rangle} \\
&=& \zeta_g\langle \pi_\nu(S_g)x_\nu, a^\ast x_\nu\rangle. 
\end{eqnarray*}
Therefore for $a\in\pi_\nu(\mathcal{O}_\Gamma)''$,
$$\langle \Delta_\nu^{1/2}\pi_\nu(S_g)x_\nu, \Delta_\nu^{1/2} a^\ast
x_\nu\rangle=\zeta_g\langle \pi_\nu(S_g)x_\nu, a^\ast x_\nu\rangle.$$
and hence for $y\in{\rm dom}(\Delta_{\nu}^{1/2})$, we have
$$\langle \Delta_\nu^{1/2}\pi_\nu(S_g)x_\nu, \Delta_\nu^{1/2} y\rangle=\zeta_g\langle \pi_\nu(S_g)x_\nu, y\rangle.$$
Thus $\Delta_\nu^{1/2}\pi_\nu(S_g)x_\nu\in{\rm
  dom}(\Delta_{\nu}^{1/2})$ and we obtain
$$\Delta_\nu\pi_\nu(S_g)x_\nu=\zeta_g\pi_\nu(S_g)x_\nu.$$
Therefore
$$\Delta_\nu^{\sqrt{-1}t}\pi_\nu(S_g)x_\nu=\zeta_g^{\sqrt{-1}t}\pi_\nu(S_g)x_\nu,$$
and then
$$(\sigma_t^\nu(\pi_\nu(S_g))-\zeta_g^{\sqrt{-1}t}\pi_\nu(S_g))x_\nu=0,$$
where $\sigma_t^\nu$ is the modular automorphism group of $\phi_\nu$.
Since $x_\nu$ is a separating vector,
$$\sigma_t^\nu(\pi_\nu(S_g))=\zeta_g^{\sqrt{-1}t}\pi_\nu(S_g).$$

Now we will show that 
$$\phi_\nu(\pi_\nu(S_g)a)=\zeta_g\phi_\nu(a \pi_\nu(S_g))
\quad\mbox{for}\quad g\in G_i\setminus H, a\in\pi_\nu(\mathcal{O}_\Gamma)''.$$
We may assume that $a=f\lambda_{g^{-1}}$ for $f\in C(\Omega)$. Recall
that $S_g=\lambda_g\chi_{\Omega\setminus Y_i}\in
C(\Omega)\rtimes_r\Gamma$. Since 
$$\phi_\nu(\pi_\nu(S_g a))=\int_{\Omega\setminus
  Y_i}f(g^{-1}\omega)d\nu(\omega)=\int_{\Omega\setminus
  Y_i}f(\omega)\frac{dg^{-1}\nu}{d\nu}(\omega)d\nu(\omega),$$
we claim that 
$$\frac{dg^{-1}\nu}{d\nu}(\omega)=\zeta_g \quad\mbox{on}\quad \Omega\setminus
  Y_i.$$ This is the Martin
  kernel $K(g^{-1}, \omega)$, (See ~\cite{woe}). Hence it suffices to show that $K(g^{-1}, x)$ is 
  constant for
  any $x=x_1\cdots x_n\in\Gamma$ such that $x_1\notin G_i$. By
  ~\cite{woe}, we have
$$K(g^{-1}, x)=\frac{G(g^{-1}, x)}{G(e, x)},$$
where $G(y, z)=\sum_{k=1}^{\infty}p^{(k)}(y, z)$ is the Green
kernel. Since any probability from $g^{-1}$ to $x$ must be through elements
of $H$ at least once, we have
$$G(g^{-1}, x)=\sum_{h\in H}F(g^{-1}, h)G(h, x),$$
where $s^x=\inf\{n\geq 0\mid 
  Z_n=x \}$ and $F(g, x)=\sum_{n=0}^{\infty}{\rm
  Pr}_g[s^x=n]$ in ~\cite{woe2}. By hypothesis $\mu(g)=\mu(hg)$
for any $g\in\bigcup_{i\in I}G_i\setminus H$ and $h\in H$, we have 
$$G(h, x)=G(e, x) \quad\mbox{for any}\quad h\in H.$$
Therefore we have $\omega_g=\log(\sum_{h\in H}F(g^{-1},
  h))$. $\sigma_t^\nu(\pi_\nu(U_h))=\pi_\nu(U_h)$ can be proved in the 
  same way. Hence we are done.  \hfill $\Box$

\medskip

\section{Appendix}

{\bf Trees \quad} We first review trees based on ~\cite{fig}. A
$graph$ is a pair $(V, E)$ consisting of a set of vertices $V$ and a family $E$ of
two-element subsets of $V$, called edges. A $path$ is a finite sequence $\{x_1,
\ldots , x_n\}\subseteq V$ such that $\{x_i, x_{i+1}\}\in E$. $(V, E)$ is said to 
be $connected$ if for $x, y\in V$ there exists a path $\{x_1,
\ldots , x_n\}$ with $x_1=x, x_n=y$. If $(V, E)$ is a tree, then for
$x, y\in V$ there 
exists a unique path $\{x_1, \ldots , x_n\}$ joining $x$ to $y$ such
that $x_i\ne x_{i+2}$. We denote this path by $[x, y]$. A tree is said to be
$locally$ $finite$ if every vertex belongs to finitely many edges. The 
number of edges to which a vertex of a locally finite tree belongs is
called a $degree$. If the degree is independent of the choice of vertices,
then the tree is called $homogeneous$. 

We introduce trees for amalgamated free product groups based on ~\cite{ser}. Let $(G_i)_{i\in I}$ be a 
family of groups with an 
index set $I$. When $H$
 is a group and every $G_i$ contains $H$ as a subgroup, then we 
 denote $\ast_H G_i$ by $\Gamma$, which is the amalgamated free
 product of the groups. If we choose sets $\Omega_i$
 of left representatives of $G_i/ H$ with $e\in\Omega_i$ for any 
$i\in I$, then each $\gamma\in\Gamma$ can be written
uniquely as
$$
   \gamma=g_1g_2 \cdots g_n h,
$$
where $h\in H, g_1\in\Omega_{i_1}\setminus\{e\}, \ldots 
,g_n\in\Omega_{i_n}\setminus\{e\}$ and $i_1 \ne i_2, i_2 \ne 
i_3, \ldots
,i_{n-1} \ne i_n$. 

Now we construct the corresponding
tree. At first, we assume that $I=\{1, 2\}$. Let
$$V=\Gamma/G_1\coprod\Gamma/G_2 \, \, \, \mbox{and}\ \, \, \, E=\Gamma/H,$$
 and the original and 
terminal maps $o :
\Gamma/H\to\Gamma/G_1$ and $t : \Gamma/H\to\Gamma/G_2$ are natural
 surjections. It is easy to see that $G_T=(V, E)$ is a tree. In general, we assume 
that the element $0$ does not belong to $I$. Let $G_0=H$ and $H_i=H$
for $i\in I$. Then we define 
$$V=\coprod_{i\in I\cup\{0\}}\Gamma/G_i \, \, \, \mbox{and}\ \, \, \, 
E=\coprod_{i\in
  I}\Gamma/H_i.$$ Now we define two maps $o, t :
E\to V$. For $H_i\in E$, let 
$$o(H_i)=G_0 \, \, \, \mbox{and}\ \, \, \, t(H_i)=G_i.$$ 
For any $\gamma H_i\in E$, we may assume that $\gamma
H=g_1\cdots g_n H_i$ such that $g_k\in\Omega_{i_k}$ with
$i_1\ne\cdots\ne i_n$. If $i=i_n$ we define
$$o(\gamma H_i)=\gamma G_{i_n} \, \, \, \mbox{and}\ \, \, \, t(\gamma
H_i)=\gamma G_0.$$
If $i\ne i_n$ we define
$$o(\gamma H_i)=\gamma G_0 \, \, \, \mbox{and}\ \, \, \, t(\gamma
H_i)=\gamma G_i.$$ 
Then we have a tree $G_T=(V, E)$. 

For a tree $(V, E)$, the set $V$ is naturally a metric space. The distance $d(x, y)$ is
defined by the number of edges in the unique path $[x, y]$. An
$infinite$ $chain$ is an infinite path $\{x_1, x_2, \ldots \}$ such that
$x_i\ne x_{i+2}$. We define an equivalence relation on the set of
infinite chains. Two infinite chains $\{x_1, x_2, \ldots\}, \{y_1, y_2,
  \ldots\}$ are equivalent if there exists an integer $k$ such that
  $x_n=y_{n+k}$ for a sufficiently large $n$. The boundary $\Omega$ of a 
  tree is the set of the equivalence classes of infinite chains. The
  boundary may be thought of as a point at infinity. Next we introduce the topology into the 
space
  $V\cup\Omega$ such that $V\cup\Omega$ is
  compact, the points of $V$ are open and $V$ is dense in
  $V\cup\Omega$. It suffices to define a basis of neighborhoods for
  each $\omega\in\Omega$. Let $x$ be a vertex. Let $\{x, x_1, x_2,
  \ldots\}$ be an infinite chain representing $\omega$. For each
  $y=x_n$, the neighborhood of $\omega$ is defined to consist of all
  vertices and all boundary points of the infinite chains which
  include $[x, y]$. 

\medskip

{\bf Hyperbolic groups \quad}We introduce hyperbolic groups defined by
Gromov. See ~\cite{gph} for details. Suppose that $(X, d)$ is a metric
space. We define a product by
$$\langle x | y \rangle_z=\frac{1}{2}\{d(x, z)+d(y, z)-d(x, y)\},$$
for $x, y, z\in X$. This is called the Gromov product. Let $\delta\geq 
0$ and $w\in X$. A metric space
$X$ is said to be $\delta$-hyperbolic with respect to $w$ if For $x, y, z\in X$,
$$\langle x | y \rangle_w\geq\min\{\langle x | z \rangle_w,  \langle y
| z \rangle_w\}-\delta. \eqno(\ddag)$$ Note that if $X$ is
$\delta$-hyperbolic with respect to $w$, then $X$ is
$2\delta$-hyperbolic with respect to any $w'\in X$. 

\begin{definition}
The space $X$ is said
to be hyperbolic if $X$ is $\delta$-hyperbolic with respect to some $w\in X$ and some
$\delta\geq 0$.  
\end{definition}

Suppose that $\Gamma$ is a
group generated by a finite subset $S$ such that $S^{-1}=S$. Let
$G(\Gamma, S)$ be the Cayley graph. The graph $G(\Gamma, S)$ has a natural word
metric. Hence $G(\Gamma, S)$ is a metric space. 

\begin{definition}
A finitely generated group $\Gamma$ is said to be hyperbolic with respect
to a finite generator system $S$ if the corresponding Cayley graph
$G(\Gamma, S)$ is hyperbolic with respect to the word metric.

In fact, hyperbolicity is independent of the choice of $S$. Therefore
we say that $\Gamma$ is a hyperbolic group, for short.
\end{definition}

We define the hyperbolic boundary of a hyperbolic space $X$. Let
$w\in X$ be a point. A sequence $(x_n)$ in $X$ is said to
$converge$ $to$ $infinity$ if $\langle x_n | x_m\rangle_w\to\infty, \, (n,
m\to\infty)$. Note that this is independent of the choice of $w$. The set $X_{\infty}$ is the 
set of all sequences converging to infinity in
$X$. Then we define an equivalence relation in $X_{\infty}$. Two
sequences $(x_n), (y_n)$ are equivalent if $\langle x_n | y_n
\rangle_w\to\infty, \, (n\to\infty)$. Although this is not an equivalence relation in 
general, the hyperbolicity assures that it is indeed an equivalence
relation. The set of all equivalent classes of $X_{\infty}$ is
called the $hyperbolic$ $boundary$ ($at$ $infinity$) and denoted by $\partial X$. Next we define
the Gromov product on $X\cup\partial X$. For $x, y\in\
X\cup\partial X$, we choose sequences $(x_n), (y_n)$ converging to $x, 
y$, respectively. Then we define $\langle x |
y\rangle=\liminf_{n\to\infty}\langle x_n | y_n\rangle_w$. Note that
this is well-defined and if $x, y\in X$ then the above product
coincides with the Gromov product on $X$. 
\begin{definition}
The topology of $X\cup\partial X$ is defined by the following neighborhood basis:
$$\{y\in X\mid d(x, y)<r\} \qquad \mbox{for}\ x\in X, r>0,$$
$$\{y\in X\cup\partial X\mid \langle x | y\rangle>r\} \qquad \mbox{for}\
x\in\partial X, r>0.$$
\end{definition}

We remark that if $X$ is a tree, then the hyperbolic boundary $\partial X$ 
coincides with the natural boundary $\Omega$ in the sense of ~\cite{fre}.

Finally we prove that an amalgamated free product $\Gamma=\ast_H
G_i$, considered in this paper, is a hyperbolic group. 

\begin{lemma}
The group $\Gamma=\ast_H G_i$ is a hyperbolic group.
\end{lemma}
{\it Proof.}
Let $S=\{g\in\bigcup_iG_i\mid |g|\leq 1\}$. Let
$G(\Gamma, S)$ be the corresponding Cayley graph. It suffices to show
$(\ddag)$ for $w=e$. For $x, y, z\in\Gamma$, we can write uniquely as
follows:
\begin{eqnarray*}
x &=& x_1\cdots x_nh_x, \\
y &=& y_1\cdots y_mh_y, \\
z &=& z_1\cdots z_kh_z,
\end{eqnarray*}
where 
$$\begin{array}{cccc}
x_1\in\Omega_{i(x_1)}, & \ldots, & x_n\in\Omega_{i(x_n)},
& h_x\in H, \\
y_1\in\Omega_{i(y_1)}, & \ldots, & y_m\in\Omega_{i(y_m)}, &
h_y\in H, \\
z_1\in\Omega_{i(z_1)}, & \ldots, & z_k\in\Omega_{i(z_k)}, &
h_z\in H.
\end{array}$$
such that each element has length one. Then $d(x, e)=n$,
$d(y, e)=m$ and $d(z, e)=k$. If $i(x_1)=i(y_1), \cdots , i(x_{l(x,
  y)})=i(y_{l(x, y)})$ and $i(x_{l(x, y)+1})\ne i(y_{l(x, y)+1})$, then 
$\langle 
x | y \rangle_e=l(x, y)$. Similarly, we obtain the positive integers
$l(x, z), l(y, x)$ such that $\langle 
x | z \rangle_e=l(x, z), \langle 
y | z \rangle_e=l(y, z)$. We can have $(\ddag)$ with $\delta=0$. \hfill $\Box$

\end{document}